\documentclass[12pt]{article}

\usepackage{amssymb,amsmath,latexsym}
\usepackage[active]{srcltx}
\usepackage[matrix,arrow,curve]{xy}
\input xy
\xyoption{all}

\oddsidemargin=\evensidemargin \advance\topmargin by-\headheight

\DeclareMathOperator{\aut}{Aut}

\DeclareMathOperator{\alt}{Alt}
\DeclareMathOperator{\bij}{Bij}

\DeclareMathOperator{\cyc}{Cyc}

\DeclareMathOperator{\Ext}{Ext}

\DeclareMathOperator{\gwr}{Gwr}

\DeclareMathOperator{\id}{id}
\DeclareMathOperator{\hol}{Hol}

\DeclareMathOperator{\inv}{Inv}

\DeclareMathOperator{\iso}{Iso}

\DeclareMathOperator{\orb}{Orb}

\DeclareMathOperator{\pr}{pr}
\DeclareMathOperator{\rad}{rad}
\DeclareMathOperator{\rk}{rk}

\DeclareMathOperator{\Span}{span}

\DeclareMathOperator{\sym}{Sym}

\def\Z{{\mathbb Z}}

\def\A{{\cal A}}
\def\B{{\cal B}}
\def\CC{{\cal C}}

\def\E{{\cal E}}

\def\H{{\cal G}}

\def\P{{\cal P}}
\def\R{{\cal R}}
\def\S{{\cal S}}
\def\scc{{\cal F}}

\def\twoe{\underset{\scriptscriptstyle ^2}{\approx}}

\def\lg{\langle}

\def\rg{\rangle}

\def\wt{\widetilde}

\def\proof{{\bf Proof}.\ }
\def\bull{\vrule height .9ex width .8ex depth -.1ex }

\makeatletter
\renewcommand{\subsection}{\@startsection{subsection}{2}{0mm}{-2mm}{-2mm}
{\bf\normalsize}}

\def\sbsnt#1{\subsection{\hspace{-3mm}#1}}
\makeatother

\newtheorem{formula}{}[section]
\newtheorem{proposition}[formula]{Proposition}
\newtheorem{definition}[formula]{Definition}
\newtheorem{corollary}[formula]{Corollary}
\newtheorem{remark}[formula]{Remark}
\newtheorem{lemma}[formula]{Lemma}
\newtheorem{theorem}[formula]{Theorem}

\newtheorem{example}[formula]{Example}

\def\thrm{\begin{theorem}}
\def\thrml#1{\begin{theorem}\label{#1}}
\def\ethrm{\end{theorem}}
\def\prpstn{\begin{proposition}}
\def\prpstnl#1{\begin{proposition}\label{#1}}
\def\eprpstn{\end{proposition}}
\def\rmrk{\begin{remark}}
\def\rmrkl#1{\begin{remark}\label{#1}}
\def\ermrk{\end{remark}}
\def\dfntn{\begin{definition}}
\def\dfntnl#1{\begin{definition}\label{#1}}
\def\edfntn{\end{definition}}
\def\nmrt{\begin{enumerate}}
\def\enmrt{\end{enumerate}}
\def\tm#1{\item[{\rm (#1)}]}
\def\qtn{\begin{equation}}
\def\qtnl#1{\begin{equation}\label{#1}}
\def\eqtn{\end{equation}}
\def\lmm{\begin{lemma}}
\def\lmml#1{\begin{lemma}\label{#1}}
\def\elmm{\end{lemma}}
\def\crllr{\begin{corollary}}
\def\crllrl#1{\begin{corollary}\label{#1}}
\def\ecrllr{\end{corollary}}
\def\hpthss{\begin{hypothesis}}
\def\hpthssl#1{\begin{hypothesis}\label{#1}}
\def\ehpthss{\end{hypotxesis}}
\def\xmpl{\begin{example}}
\def\xmpll#1{\begin{example}\label{#1}}
\def\exmpl{\end{example}}
\def\css{\begin{cases}}
\def\ecss{\end{cases}}

\begin{document}

\title{Schurity of S-rings over a cyclic group and\\
generalized wreath product of\\ permutation groups}
\author{
Sergei Evdokimov \\[-2pt]
\small Steklov Institute of Mathematics\\[-4pt]
\small at St. Petersburg \\[-4pt]
{\tt \small evdokim@pdmi.ras.ru }
\thanks{The work was partially supported by
Slovenian-Russian bilateral project, grant no. BI-RU/10-11-018.}
\and
Ilya Ponomarenko\\[-2pt]
\small Steklov Institute of Mathematics\\[-4pt]
\small at St. Petersburg \\[-4pt]
{\tt \small inp@pdmi.ras.ru}
\thanks{The work was partially supported by
Slovenian-Russian bilateral project, grant no. BI-RU/10-11-018.
}
}
\date{}

\maketitle

\begin{abstract}
The generalized wreath product of permutation groups is introduced. By means of it
we study the schurity problem for S-rings over a cyclic group $G$ and the automorphism groups
of them. Criteria for the schurity and non-schurity of the generalized wreath product of two
such S-rings are obtained. As a byproduct of the developed theory we prove that
$G$ is a Schur group whenever the total number $\Omega(n)$ of prime factors of the
integer~$n=|G|$ is at most~$3$. Moreover, we describe the structure
of a non-schurian S-ring over $G$ when $\Omega(n)=4$. The latter result implies in particular
that if $n=p^3q$ where $p$ and $q$ are primes, then $G$ is a Schur group.
\end{abstract}

\section{Introduction}\label{240111b}
One of the central problems in the theory of S-rings over a finite group that goes back to
H.~Wielandt, is to identify schurian S-rings, i.e. those arising from suitable permutation
groups (as for a background of S-rings see Section~\ref{081209i}). At present this problem
is still open even for S-rings over a cyclic group which are called below {\it circulant}.
In \cite{LM96,EP01ce} it was proved
that any circulant S-ring can be constructed from S-rings of rank~$2$ and normal S-rings
by means of two operations: tensor product and generalized wreath product. It should be noted
that the S-rings of rank~$2$ as well as the normal S-rings are schurian, and the tensor product
of schurian S-rings is schurian. However, the generalized wreath product of schurian S-rings is
not always schurian: the first examples of non-schurian generalized wreath products
were constructed in~\cite{EP01ae}. It follows that there exist non-schurian S-rings
over a cyclic group of order~$n$ where $\Omega(n)=4$.  On the
other hand, it was proved in \cite{P74,KP} that any cyclic group of order~$n$ with
$\Omega(n)\le 2$ is a {\it Schur group}, i.e. any S-ring over it is schurian.\medskip

The main goal of this paper is to lay the foundation of the complete characterization of
Schur cyclic groups. As an easy byproduct of the theory developed here we obtain
the following result (which is just a reformulation of Theorem~\ref{140510a}).

\thrml{270910c}
Any cyclic group of order $n$ with $\Omega(n)\le 3$ is a Schur group.
\ethrm

The problem of the identification of schurian S-rings leads naturally to studying
the automorphism group $\aut(\A)$ of an S-ring $\A$ over a group~$G$. By definition
$\aut(\A)$ is the automorphism group of the Cayley scheme corresponding to~$\A$; this
group always contains the subgroup $G_{right}$ consisting of all permutations induced by right multiplications.
If the group~$G$ is cyclic, then the arguments from the above paragraph show that
all we have to do is to study $\aut(\A)$ when $\A$ is a
generalized wreath product:
\qtnl{280311a}
\A=\A_U\wr_{U/L}\A_{G/L}
\eqtn
where $U$ and $L$ are $\A$-subgroups of the group $G$ such that $L\le U$, and
$\A_U$ (resp. $\A_{G/L}$) is the restriction of~$\A$ to $U$ (resp. to~$G/L$);
the latter is equivalent to the fact that any basic set of~$\A$ contained in
$G\setminus U$ is a union of cosets modulo~$L$. In these terms the group $\aut(\A)$
coincides with the largest group $\Gamma\le\sym(G)$ such that
$$
\Gamma^U\le\aut(\A_U)\quad\text{and}\quad\Gamma^{G/L}\le\aut(\A_{G/L})
$$
where $\Gamma^U$ (resp. $\Gamma^{G/L}$) is the permutation group induced by the 
action on $U$ of
the setwise stabilizer of $U$ in $\Gamma$ (resp. by the action on $G/L$ of the group $\Gamma$).
This suggests the following construction which is central for our paper.\medskip

Let $V$ be a finite set and $E_0\subset E_1$ equivalence relations on~$V$. Suppose we
are given two permutations groups $\Delta_U\le\sym(U)$ where $U\in V/E_1$, and $\Delta_0\le\sym(V/E_0$,
such that
$$
(\Delta_U)^{U/E_0}=(\Delta_0)^{U/E_0}
$$
where $U/E_0$ is the set of $E_0$-classes contained in~$U$. Then all maximal groups
$\Gamma\le\sym(V)$ leaving $E_0$ and $E_1$ fixed and such that
$$
\Gamma^U=\Delta_U\quad\text{and}\quad\Gamma^{V/E_0}=\Delta_0,
$$
are permutationally isomorphic. In fact,
any such group is uniquely determined by choosing a family of bijections
identifying the classes of $V/E_1$ with~$U$. It is called the {\it generalized wreath
product} of $\Delta_U$ and $\Delta_0$ (with respect to this family). It is easily seen
that our construction coincides with the ordinary wreath product when $E_0=E_1$.\medskip

An important special case arises when $V=G$ is a group, the classes of~$E_0$ and $E_1$ are
the left cosets modulo a normal subgroup $L$ of~$G$ and modulo a subgroup $U$ of~$G$ that
contains~$L$ respectively. Then one can identify the classes of $V/E_1$ with~$U$ by means of
permutations from $G_{right}$, and the generalized generalized wreath product of $\Delta_U$
and $\Delta_0$ does not depend on choosing bijections of this type. The permutation
group on $G$ obtained in this way contains~$G_{right}$, is called the {\it canonical} generalized
wreath product of~$\Delta_U$ by~$\Delta_0$ over~$G$, and is denoted by
$\Delta_U\wr_{U/L}\Delta_0$ (see Subsection~\ref{130111a}).\medskip

In the early 1950's D.~K.~Faddeev introduced the generalized wreath
product of two abstract groups (see \cite[p.46]{ILF}). One can prove that the generalized
wreath product of permutation groups defined in the above paragraph is isomorphic
(as an abstract group) to the group obtained by the Faddeev construction. Curiously,
the authors independently came to the generalized wreath product of S-rings and only after that
got aware of the Faddeev work. It should be noted  that the concept of generalized
wreath product of permutation groups proposed in paper \cite{BPRS} is different from ours.\medskip

One of the motivations to define the generalized wreath product of permutation groups
is Theorem~\ref{270910d} below which immediately follows from Corollaries~\ref{050710a}
and~\ref{260210a}. In particular, this gives a necessary and sufficient condition for a generalized
wreath product of two S-rings to be schurian in terms of $2$-equivalent subgroups
of their automorphism groups (for the $2$-equivalence we use the same notation $\twoe$
as in Wielandt's book~\cite{W69}, see also Notation).

\thrml{270910d}
Let $\A$ be an S-ring over an abelian group $G$ that is the generalized wreath
product~\eqref{280311a}. Then 
\nmrt
\tm{1} $\aut(\A)=\aut(\A)^U\wr_{U/L}\aut(\A)^{G/L}$,
\tm{2} $\A$ is schurian if and only if there exist groups $\Delta_U$ and $\Delta_0$ such that
\nmrt
\tm{a} $U_{right}\le\Delta_U\le\aut(\A_U)$ and $(G/L)_{right}\le\Delta_0\le\aut(\A_{G/L})$,
\tm{b} $\Delta_U\twoe\aut(\A_U)$ and $\Delta_0\twoe\aut(\A_{G/L})$,
\tm{c} $(\Delta_0)^{U/L}=(\Delta_U)^S$.
\enmrt
Moreover, in this case $\aut(\A)\twoe\Delta_U\wr_{U/L}\Delta_0$.
\enmrt
\ethrm

The proof of Theorem~\ref{270910d} is not too complicated modulo the characterization
of the automorphisms of the generalized wreath products of S-rings given in paper
\cite{EP03be}. 
However, the theory of generalized wreath products of  circulant S-rings
initiated in that paper and developed here
(Sections~\ref{041010a}-\ref{041010d}), especially the part of it concerning
singularities and their resolutions, enables us to reveal a special case of
Theorem~\ref{270910d} in which the schurity is easily checked. This special
case formulated in the theorem below is, in fact, the key ingredient to
prove Theorem~\ref{270910c} and the other results of this paper.

\thrml{190210a}
Let $\A$ be an S-ring over a cyclic group~$G$ that is the generalized wreath
product~\eqref{280311a}. Suppose that $\A_{U/L}$ is the tensor product of a normal
S-ring and S-rings of rank~2. Then $\A$ is schurian if and only if so are the S-rings~$\A_U$
and~$\A_{G/L}$.
\ethrm

According to Theorem~\ref{mbnulc} of this paper first proved in~\cite{EP01ce}, any circulant S-ring
with trivial radical is the tensor product of a normal
S-ring and S-rings of rank~2. Thus the following statement is a special case of
 Theorem~\ref{190210a}.

\crllrl{230810c}
Let $\A$ be an S-ring over a cyclic group~$G$ that is the generalized wreath
product~\eqref{280311a}. Suppose
that $\A_{U/L}$ is an S-ring with trivial radical. Then $\A$ is schurian if and only if so
are the S-rings~$\A_U$ and $\A_{G/L}$.
\ecrllr

With the exception of the Burnside-Schur theorem it is not so much known on the structure of a
permutation group containing a regular cyclic subgroup. In fact, only the primitive case was
studied in detail and the main results are based on the classification of finite simple groups.
On the other hand, as a byproduct of the theory developed to prove Theorem~\ref{190210a} we can
get some information on a $2$-closed permutation group containing a regular cyclic subgroup
(such a group is nothing else than the automorphism group of a circulant S-ring). Namely, in
Section~\ref{051010a} we prove the following result.

\thrml{240810a}
Let $\Gamma$ be the automorphism group of an S-ring over a cyclic group~$G$. Then
\nmrt
\tm{1} every non-abelian composition factor of~$\Gamma$ is an alternating group,
\tm{2} the group $\Gamma$ is $2$-equivalent to a solvable group containing $G_{right}$
if and only if every alternating composition factor of~$\Gamma$ is of prime degree.
\enmrt
\ethrm

By Theorem~\ref{270910c} non-schurian S-rings over a cyclic group of order~$n$ exist
only if $\Omega(n)\ge 4$. In Section~\ref{051010b} we study the structure of a non-schurian
S-ring $\A$ when
$\Omega(n)=4$. Based on Theorem~\ref{190210a} we prove that $\A$ is the generalized
wreath product of two circulant S-rings such that either both of them are
generalized wreath products, or exactly one of them is normal. We note that non-schurian
circulant S-rings constructed in~\cite{EP01ae} are of the former type. In Section~\ref{051010g}
non-schurian circulant S-rings of the latter type are presented. Finally, Theorem~\ref{211210a}
and statement~(1) of Theorem~\ref{051010c} imply the following result.

\thrml{230810b}
Any cyclic group of order $p^3q$ where $p$ and $q$ are primes, is a Schur group.\bull
\ethrm

Concerning permutation groups we refer to~\cite{W64} and~\cite{W69}.
For the reader convenience we collect the basic facts on S-rings and Cayley schemes over an abelian
group in Section~\ref{011209b}. The theory of S-rings over a cyclic group developed by
the authors in~\cite{EP01ae,EP01ce,EP03be} is given in Section~\ref{071010a}.
\medskip

{\bf Notation.}
As usual by $\Z$ we denote the ring of integers.\medskip

For a positive integer $n$ we denote by $\Omega(n)$  the total number of prime factors of~$n$.\medskip

The set of all equivalence relations on a set $V$ is denoted by $\E(V)$.\medskip

For $X\subset V$ and $E\in\E(V)$ we set $X/E=X/E_X$ where $E_X=X^2\cap E$.
If the classes of $E_X$ are singletons, we identify $X/E$ with $X$.\medskip

For $R\subset V^2$, $X\subset V$ and $E\in\E(V)$ we set
$$
R_{X/E}=\{(Y,Z)\in(X/E)^2: R_{Y,Z}\ne\emptyset\}
$$
where $R_{Y,Z}=R\cap (Y\times Z)$.\medskip

The set of all bijections from $V$ onto $V'$ is denoted by $\bij(V,V')$.\medskip

For $B\subset\bij(V,V')$, $X\subset V$ and $X'\subset V'$, $E\in\E(V)$ and $E'\in\E(V')$ we set
$$
B^{X/E,X'/E'}=\{f^{X/E}:\ f\in B,\ X^f=X',\ E^f=E'\},
$$
$$
B^{X/E}=B^{X/E,X/E}
$$
where $f^{X/E}$ is the bijection from $X/E$ onto $X'/E'$ induced by~$f$.\medskip

The group of all permutations of $V$ is denoted by $\sym(V)$.\medskip

The set of orbits of a group $\Gamma\le\sym(V)$ is denoted by $\orb(\Gamma)=\orb(\Gamma,V)$.\medskip

The setwise stabilizer of a set $U\subset V$ in the group~$\Gamma$ is denoted by $\Gamma_{\{U\}}$.\medskip

For an equivalence relation $E\in\E(V)$ we set $\Gamma_E=\bigcap_{X\in V/E}\Gamma_{\{X\}}$.\medskip

We write $\Gamma\twoe\Gamma'$ if groups $\Gamma,\Gamma'\le\sym(V)$ are $2$-equivalent, i.e.
have the same orbits in the coordinatewise action on~$V^2$.\medskip

The permutation group defined by the right multiplications of a group~$G$ on itself is
denoted by $G_{right}$.\medskip

The holomorph $\hol(G)$ is identified with the permutation
group on the set~$G$ generated by $G_{right}$ and~$\aut(G)$.\medskip

The natural epimorphism from $G$ onto $G/H$ is denoted by $\pi_{G/H}$.\medskip

For sections $S=L_1/L_0$ and $T=U_1/U_0$ of $G$ we write $S\le T$, if
$U_0\le L_0$ and $L_1\le U_1$.

\section{Schemes and S-rings}\label{011209b}

\sbsnt{Schemes.}\label{250706e}
As for the background on scheme theory presented here see~\cite{EP09} and
references there.
Let $V$ be a finite set and $\R$ a partition of $V^2$.
Denote by $\R^*$ the set of all unions of the elements of $\R$.\medskip

A pair $\CC=(V,\R)$ is
called a {\it coherent configuration} or a {\it scheme} on $V$ if the following
conditions are satisfied:
\nmrt
\item[$\bullet$] the diagonal of $V^2$ belongs to $\R^*$,
\item[$\bullet$] $\R$ is closed with respect to transposition,
\item[$\bullet$] given  $R,S,T\in\R$, the number $|\{v\in V:\,(u,v)\in R,\ (v,w)\in S\}|$
does not depend on the choice of $(u,w)\in T$.
\enmrt

The elements of $V$, $\R=\R(\CC)$, $\R^*=\R^*(\CC)$ and the number from the third condition,
are called the {\it points}, the {\it basis relations}, the {\it relations} and
the {\it intersection number} (associated with $R,S,T$) of~$\CC$, respectively. The
numbers $|V|$ and $|\R|$ are called the {\it degree} and {\it rank} of~$\CC$.
If the diagonal of $V^2$ belongs to~$\R$, the scheme $\CC$ is called {\it homogeneous}.\medskip

Two schemes are called {\it isomorphic} if there exists a bijection between their point
sets preserving the basis relations. Any such bijection is called an {\it isomorphism}
of these schemes. The group of all isomorphisms of a scheme $\CC$ contains a normal
subgroup
$$
\aut(\CC)=\{f\in\sym(V):\ R^f=R,\ R\in\R\}
$$
called the {\it automorphism group} of~$\CC$. If $V$ coincides with a group $G$ and
$G_{right}\le\aut(\CC)$, then $\CC$ is called a {\it Cayley scheme} over~$G$.
Such a scheme is {\it normal} if the group $G_{right}$ is normal in~$\aut(\CC)$.\medskip

Given a permutation group $\Gamma\le\sym(V)$ set $\orb_2(\Gamma)=\orb(\Gamma,V^2)$
where the latter is the set of orbits in the coordinatewise action of $\Gamma$ on~$V^2$.
Then the pair
$$
\inv(\Gamma)=(V,\orb_2(\Gamma))
$$
is a scheme; we call it the {\it scheme of the group} $\Gamma$. Any scheme of this
type is called {\it schurian}. One can see that
$\Gamma\le\aut(\CC)$ where где $\CC=\inv(\Gamma)$.\medskip

Let $\CC$ be a homogeneous scheme and $\E(\CC)=\R^*(\CC)\cap\E(V)$. Then given
a class $X$ of an equivalence relation belonging to $\E(\CC)$ and $E\in\E(\CC)$, the pair
$\CC_{X/E}=(X/E,\R_{X/E})$
where $\R_{X/E}=\{R_{X/E}:\ R\in\R,\ R_{X/E}\ne\emptyset\}$,
is a scheme on $X/E$. One can see that
$$
\aut(\CC)^{X/E}\le\aut(\CC_{X/E}).
$$
Moreover, if the scheme $\CC$ is schurian, then the above two groups are $2$-equivalent.
When $\CC$ is a Cayley scheme over a group $G$ and $V/E=G/L$ for
a normal subgroup $L$ in~$G$, we set $\CC_{X/L}=\CC_{X/E}$.\medskip

Schemes $\CC$ and $\CC'$ are called {\it similar} if there exists a bijection
$\varphi:\R\to\R'$ called a {\it similarity} from~$\CC$ to~$\CC'$,
such that the intersection number associated with $R,S,T\in\R$ equals
the intersection number associated with $R^\varphi,S^\varphi,T^\varphi$. The set
of all isomorphisms from~$\CC$ to~$\CC'$ inducing
a similarity~$\varphi$ is denoted by $\iso(\CC,\CC',\varphi)$. Let $\CC$ be a
homogeneous scheme and $E\in\E(\CC)$. Then given $X,X'\in V/E$ there
exists a uniquely determined bijection
\qtnl{300810l}
\varphi_{X,X'}:\R_{X^{}}\to\R_{X'},\ R_{X^{}}\mapsto R_{X'}.
\eqtn
Moreover, this bijection defines a similarity from $\CC_{X^{}}$ to $\CC_{X'}$.

\sbsnt{Schur rings and Cayley schemes.}\label{081209i}
Let $G$ be a finite group. A subring~$\A$ of the group ring~$\Z G$ is
called a {\it Schur ring} ({\it S-ring}, for short) over~$G$ if it has a
(uniquely determined) $\Z$-basis consisting of the elements
$\sum_{x\in X}x$ where $X$ runs over the classes of a partition $\S=\S(\A)$ of~$G$ such that
$$
\{1\}\in\S,\quad
\textstyle{\rm and}\quad
X\in\S\ \Rightarrow\ X^{-1}\in\S.
$$
Let $\A'$ be an S-ring over a group $G'$. A group isomorphism $f:G\to G'$ is called a
{\it Cayley isomorphism} from $\A$ to $\A'$ if $\S(\A)^f=\S(\A')$.\medskip

The elements of $\S$ and the number $\rk(\A)=|\S|$ are called respectively the {\it basic}
sets and the {\it rank} of the S-ring~$\A$.
Any union of basic sets is called an
{\it $\A$-subset of~$G$} or {\it $\A$-set}; the set of all of them is denoted by $\S^*(\A)$.
It is easily seen that this set is closed with respect to taking inverse and product. Given
$X\in\S^*(\A)$ we set
$$
\S(\A)_X=\{Y\in\S(\A):\ Y\subset X\}.
$$
The $\Z$-submodule of~$\A$ spanned by this set is denoted by $\A_X$.
If $\A'$ is an S-ring over~$G$ such that $\S^*(\A)\subset\S^*(\A')$, then we write $\A\le\A'$.\medskip

Any subgroup of the group $G$ that is an $\A$-set is called an {\it $\A$-subgroup of~$G$} or
{\it $\A$-group}; the set of all of them is denoted by $\H(\A)$. The S-ring $\A$ is called
{\it dense} if every subgroup of~$G$ is an $\A$-group, and {\it primitive} if the only
$\A$-subgroups are $1$ and $G$.\medskip

If $\A_1$ and $\A_2$ are S-rings over groups $G_1$ and $G_2$ respectively, then
the subring $\A=\A_1\otimes \A_2$ of the ring $\Z G_1\otimes\Z G_2=\Z G$ where
$G=G_1\times G_2$, is an S-ring over the group $G$ with
$$
\S(\A)=\{X_1\times X_2: X_1\in\S(\A_1),\ X_2\in\S(\A_2)\}.
$$
It is called the {\it tensor product} of $\A_1$ and $\A_2$. The following
statement was proved in~\cite{EP10}.

\lmml{130209d}
Let $\A$ be an S-ring over the direct product $G$ of groups $G_1$ and $G_2$ such that
$G_1,G_2\in\H(\A)$. Then $\pi_i(X)\in\S(\A)$ for all $X\in\S(\A)$
where $\pi_i$ is the projection of $G$ on $G_i$, $i=1,2$. In particular,
$\A\ge\A_{G_1}\otimes \A_{G_2}$.\bull
\elmm

Given a group~$G$ there is a 1-1 correspondence between the S-rings over~$G$ and the Cayley
schemes over~$G$ that preserves the natural partial orders on these sets: any basis
relation of the scheme $\CC$ corresponding to an S-ring $\A$ is of the form
$$
\{(g,xg):\ g\in G,\ x\in X\},\qquad X\in\S(\A).
$$
It can be proved that given a group $H\le G$ we have $H\in\H(\A)$ if and only if $E_H\in\E(\CC)$
where $E_H$ is the equivalence relation on $G$ with
$$
G/E_H=\{Hg:\ g\in G\}.
$$
The group $\aut(\A):=\aut(\CC)$ is called the {\it automorphism group} of the S-ring~$\A$.
An S-ring is {\it schurian}
(resp. {\it normal}) if so is the corresponding Cayley scheme.

\section{Sections in groups and S-rings}\label{041010a}

\sbsnt{Sections of a group.} Let $G$ be a group. Denote by $\scc(G)$ the set of all its sections,
i.e. quotients of subgroups of~$G$. A section $U_1/U_0$ is called a {\it multiple} of a section
$L_1/L_0$ if
$$
L_0=U_0\cap L_1\quad\text{and}\quad U_1=U_0L_1.
$$
Let us define an equivalence relation on the set $\scc(G)$
as the transitive closure of the relation "to be a multiple". The set of all
equivalence classes is denoted by $\P(G)$. Any two elements of the same equivalence class
are called {\it projectively equivalent}. Trivial $G$-sections, i.e. sections of
order~$1$, are obviously projectively equivalent.\medskip

Let $S=L_1/L_0$ and $T=U_1/U_0$ be projectively equivalent $G$-sections.
If $T$ is a multiple of $S$, then obviously the groups $S=L_1/U_0\cap L_1$
and $T=U_0L_1/U_0$ are isomorphic under the natural isomorphism
\qtnl{251010a}
f_{S,T}:S\to T,\quad gL_0\mapsto gU_0.
\eqtn
Generally, the sections $S$ and $T$ remain isomorphic and an isomorphism can be defined
by a suitable composition of the above isomorphisms and their inverses. Any such isomorphism
will be called a {\it projective} isomorphism from $S$ onto~$T$.\medskip

Let $\Gamma\le\sym(G)$ be a group containing $G_{right}$. A section $U/L\in\scc(G)$ is called
{\it $\Gamma$-invariant} if $E_U$ and $E_L$ are $\Gamma$-invariant equivalence relations. In this
case we write $\Gamma^{U/L}$ instead of $\Gamma^{U/E_L}$, and $\gamma^{U/L}$ instead of
$\gamma^{U/L}$, $\gamma\in\Gamma$. The following easy lemma is very
useful.

\lmml{160210a}
Let $G$ be a group and $G_{right}\le\Gamma\le\sym(G)$.
Then given projectively equivalent $\Gamma$-invariant $G$-sections $S$ and $T$,
any projective isomorphism $f:S\to T$ induces a permutation isomorphism from $\Gamma^S$ onto
$\Gamma^T$. Moreover,
$$
(S_{right})^f=T_{right}
$$
and $(\gamma^S)^f=\gamma^T$ for all $\gamma\in\Gamma$ leaving the point~$1_G$ fixed.
\elmm
\proof Without loss of generality we can assume that $T=U_1/U_0$ is a multiple
of $S=L_1/L_0$ and $f=f_{S,T}$ where $f_{S,T}$ is defined in~\eqref{251010a}.
Then given a permutation $\gamma\in\Gamma_{\{L_1\}}$ and $g\in L_1$
the block $(gL_0)^\gamma$ contains the element $g^\gamma$. It follows that
$(gL_0)^\gamma=g^\gamma L_0$, and hence by the definition of $f$ we have
$$
((gL_0)^\gamma)^f=(g^\gamma L_0)^f=g^\gamma U_0.
$$
Thus $(\Gamma^S)^f\le\Gamma^T$ and $S_{right}^f=T_{right}$. Moreover, this also proves
the statement for $\gamma$ leaving the point~$1_G$ fixed.
Conversely, let $\gamma\in\Gamma_{\{U_1\}}$.
Then by the above it suffices to find $\gamma'\in\Gamma_{\{U_1\}}$ such that $(\gamma')^T=\gamma^T$
and $\gamma'\in\Gamma_{\{L_1\}}$. However, since the equivalence relation $E_{L_1}$ is
$\Gamma$-invariant, we have $L_1^\gamma=L_1g$ for some $g\in U_1$. Taking into account
that $U_1=L_1U_0$, we can assume that $g\in U_0$. Set $\gamma'=\gamma\gamma_1$
where $\gamma_1$ is a permutation on $G$ taking $x$ to $xg^{-1}$. Then
$$
(L_1)^{\gamma'}=L_1,\qquad (\gamma')^T=\gamma^T.
$$
Since obviously $\gamma'\in\Gamma_{\{U_1\}}$, we are done.\bull

\sbsnt{Sections of an S-ring.}
Let $\A$ be an S-ring over a group~$G$. Denote by $\scc(\A)$ the set of all {\it $\A$-sections},
i.e. those $G$-sections $S=U/L$ for which both $U$ and $L$ are $\A$-groups.
Given such an $\A$-section and a set $X\in\S(\A)_U$ put $X_S=\pi_{U/L}(X)$. Then the $\Z$-module
$$
\A_S=\Span \{X_S:\ X\in\S(\A)_U\}
$$
is an S-ring over the group~$S$ the basic sets of which are exactly the sets under the above
span. We say that the section $S$ is of rank~$2$ (resp. normal, primitive) if so is the S-ring
$\A_S$.\medskip

Denote by $\P(\A)$ the set of nonempty sets $C\cap\scc(\A)$ where $C\in\P(G)$. Then $\P(\A)$
forms a partition of the set $\scc(\A)$ into classes  of projectively equivalent $\A$-sections.
 It should be noted that if $\A'\ge\A$, then
$\scc(\A')\supset\scc(\A)$ and each class of projectively equivalent $\A$-sections is contained
in a unique class of projectively equivalent $\A'$-sections.

\thrml{261010a}
Let $\A$ be an S-ring over a group $G$. Then given projectively equivalent $\A$-sections
$S=L_1/L_0$ and $T=U_1/U_0$ any projective isomorphism $f:S\to T$ is a Cayley isomorphism
from $\A_S$ onto $\A_T$. Moreover, if $T$ is a multiple of $S$, then
$$
(X_S)^f=X_T^{},\qquad X\in\S(\A)_{L_1}.
$$
\ethrm
\proof Without loss of generality we can assume that~$T$ is a multiple of~$S$,
and $L_0=1$, $U_1=G$. Then the projective isomorphism $f_{S,T}$ from~\eqref{251010a}
coincides with the bijection defined by formula~(6) in~\cite{EP03be}
for $V=G$, $E_1=E_{L_1}$, $X_1=L_1$ and $E_2=E_{U_0}$. Thus the
required statement follows from statement~(2) of Theorem~2.2 of that paper
for $\CC$ being the Cayley scheme associated with~$\A$.\footnote{In fact, Theorem~2.2 was
proved in~\cite{EP03be} for a commutative scheme. However, statement~(2) of it is true under a weaker
assumption that the adjacency matrices of the equivalence relations
$E_1$ and $E_2$ commute. The latter is true because $L_1U_0=U_0L_1$.}\bull\medskip

Obviously, any two sections of a class $C\in\P(G)$ have the same order; we call it the {\it order}
of~$C$. If, in addition, $C\in\P(\A)$, then from Theorem~\ref{261010a} it follows that all 
sections in~$C$ have the same rank~$r$, and if $C$ contains a primitive (resp. normal) 
section, then all sections in~$C$ are primitive (resp. normal). In these cases we say that 
the class $C$ is a class of rank~$r$, or a primitive (resp. normal) class.\medskip

\crllrl{160210b}
Let $\A$ be an  S-ring over a group~$G$ and $C$ a class
of projectively equivalent $\A$-sections. Suppose that
there exists a section $S\in C$ such that
\qtnl{190410a}
\aut(\A)^S\le\hol(S).
\eqtn
Then this inclusion holds for all $S\in C$.
\ecrllr
\proof Let $T\in C$. Then by Lemma~\ref{160210a} for $\Gamma=\aut(\A)$ there exists
a permutation isomorphism of the group $\aut(\A)^S$ onto the group $\aut(\A)^T$.
This isomorphism takes $S_{right}$ to $T_{right}$, and hence $\hol(S)$ to
$\hol(T)$. Thus $\aut(\A)^T\le\hol(T)$.\bull\medskip

\sbsnt{The $S$-condition.}
Let $\A$ be an S-ring over a group $G$ and $S=U/L$ an $\A$-section. Following~\cite{EP01ce} we
say that $\A$ satisfies the {\it $S$-condition}, if $L$ is
normal in~$G$ and
$$
LX=XL=X,\qquad X\in\S(\A)_{G\setminus U}.
$$
If, moreover, $L\ne 1$ and $U\ne G$, we say that $\A$ satisfies the $S$-condition
{\it nontrivially}. The following theorem is a specialization
of \cite[Theorem~3.1]{EP01ae}.

\thrml{160710a}
Let $S=U/L$ be a section of an abelian group~$G$, and let $\A_1$ and $\A_2$ be S-rings
over groups $U$ and $G/L$ respectively. Suppose that
$$
S\in\scc(\A_1)\cap\scc(\A_2)\quad\text{and}\quad(\A_1)_S=(\A_2)_S.
$$
Then there is a uniquely
determined S-ring~$\A$ over the group~$G$ that satisfies the $S$-condition and such that
$\A_U=\A_1$ and $\A_{G/L}=\A_2$.\bull
\ethrm

The S-ring $\A$ from Theorem~\ref{160710a} is called the {\it generalized wreath product}
of S-rings $\A_1$ and $\A_2$; we denote it by $\A_1\wr_S\A_2$ and omit $S$ and the
word "generalized"\  when $|S|=1$. Thus an S-ring $\A$ over~$G$
satisfies the $S$-condition if and only if $\A=\A_U\wr_S\A_{G/L}$. We say that $\A$ is a
{\it nontrivial} or {\it proper} generalized wreath product if it satisfies the $S$-condition
nontrivially.

\section{S-rings over a cyclic group}\label{071010a}

\sbsnt{General theory.}
Let $\A$ be an S-ring over a cyclic group~$G$. One can see that given $X\subset G$ the set
$$
\rad(X)=\{g\in G:\ gX=Xg=X\}
$$
is an $\A$-subgroup of~$G$ whenever $X\in\S^*(\A)$. From the well-known Schur theorem
on multipliers \cite[Theorem~23.9]{W64} it follows that the group $\rad(X)$
does not depend on the choice of $X\in\S(\A)$ such that $X$ contains a generator
of~$G$. We call this group the {\it radical} of~$\A$ and denote it by $\rad(\A)$.

\thrml{mbnulc}
Let $\A$  be a circulant S-ring. Then
\nmrt
\tm{1} $\rad(\A)\ne 1$ if and only if $\A$ is a nontrivial generalized wreath product,
\tm{2} $\rad(\A)= 1$ if and only if $\A$ is the tensor product of a normal S-ring with trivial
radical and S-rings of rank~$2$.
\enmrt
\ethrm
\proof Follows from \cite[Corollaries~5.5,6.4]{EP01ce}.\bull\medskip

The S-ring $\A$ is called {\it cyclotomic}\footnote{In~\cite{EP01ce} such an S-ring was
called an orbit one.} if $\S(\A)=\orb(K,G)$ for some group $K\le\aut(G)$; in this case we
write $\A=\cyc(K,G)$.

\thrml{030910a}
Let $\A$  be a normal S-ring over a cyclic group~$G$. Then
\nmrt
\tm{1} $\A$ is cyclotomic; in particular, $\A$ is dense and schurian,
\tm{2} $|\rad(\A)|=1$ or $2$,
\tm{3} if $\rad(\A)\ne 1$, then $\A=\A_U\wr_{U/L}\A_{G/L}$ where $|G/U|=|L|=2$ and $\A_U$ and $\A_{G/L}$
are normal S-rings with trivial radicals.
\enmrt
\ethrm
\proof Follows from \cite[Theorems~6.1,5.7]{EP01ce}.\bull

\crllrl{170510b}
Any circulant S-ring with trivial radical is schurian.
\ecrllr
\proof By statement~(2) of Theorem~\ref{mbnulc} it suffices to verify that any normal
circular ring is schurian. However, this is true by statement~(1) of Theorem~\ref{030910a}.\bull\medskip

It is well known that
\qtnl{060910a}
\aut(G)=\prod_{p\in\P}\aut(G_p)
\eqtn
where $\P$ is the set of primes dividing $|G|$ and $G_p$ is the Sylow $p$-subgroup of~$G$. Given
$p\in\P$ the group $\aut(G_p)$ is a cyclic group of order $(p-1)|G_p|/p$ whenever $p$ is odd.
In this case it is easily seen that if $X$ is a orbit of a group $K\le\aut(G_p)$ containing a
generator of~$G_p$, then $\rad(X)\ne 1$ if and only if $p$ divides $|K|$. Thus the cyclotomic
S-ring $\cyc(K,G_p)$ has trivial radical if and only if $|K|$ is coprime to~$p$.




\crllrl{091210v}
Let $\A$ be a normal S-ring over a cyclic $p$-group and $S$ an $\A$-section. Suppose
that $p$ is odd and $|S|\ge p^2$. Then the S-ring $\A_S$ is normal.
\ecrllr
\proof Since $p>2$, statements~(1) and (2) of Theorem~\ref{030910a} imply that
$\A=\cyc(K,G)$ where $K\le\aut(G)$, and $\rad(\A)=1$. In particular, $p$ does not
divide $|K|$. Therefore $\A_S=\cyc(K^S,G)$ and $\rad(\A_S)=1$ (because $|K^S|$ is coprime to~$p$).
Since $\rk(\A_S)>2$, this implies by statement~(2) of Theorem~\ref{mbnulc} that the S-ring $\A_S$
is normal.\bull\medskip

\sbsnt{Singularities.}\label{170111a}
According to \cite{EP03be} we say that an $\A$-section $S$ is the smallest (respectively, 
greatest) if every section which is pojectively equivalent to~$S$ is a multiple of~$S$
(respectively, if $S$ is a multiple of every section which is pojectively equivalent to~$S$).
The following result was proved in Lemma~5.2 of that paper.

\thrml{281010a}
Let $\A$  be a circulant S-ring. Then any class of projectively equivalent
$\A$-sections contains the smallest and greatest elements.\bull
\ethrm

A class of projectively equivalent $\A$-sections of rank~$2$ is called {\it singular} if
its order is more than~$2$ and it contains two sections $S=L_1/L_0$ and $T=U_1/U_0$ such that $T$ is a multiple of~$S$ and
the following two conditions are satisfied:
\nmrt
\tm{S1} $\A=\A_{U_0}\wr_{U_0/L_0}\A_{G/L_0}=\A_{U_1}\wr_{U_1/L_1}\A_{G/L_1}$,
\tm{S2} $\A_{U_1/L_0}=\A_{L_1/L_0}\otimes\A_{U_0/L_0}$.
\enmrt

\thrml{160210c}
Let $\A$ be an S-ring over a cyclic group~$G$. Then for any primitive class
$C\in\P(\A)$ one of the following statements holds:
\nmrt
\tm{1} $C$ is singular and $\aut(\A)^S=\sym(S)$ for all $S\in C$,
\tm{2} $C$ is of prime order and $\aut(\A)^S\le\hol(S)$ for all $S\in C$.
\enmrt
\ethrm
\proof First, suppose that the class $C$ contains a {\it subnormal} section~$S$; by
definition this means that $S$ is (in the natural sense) a section of a
normal $\A$-section. Then $\aut(\A)^S\le\hol(S)$. By Corollary~\ref{160210b} this implies that
$\aut(\A)^{S'}\le\hol(S')$ for all $S'\in C$. Moreover, by statement~(1) of Theorem~\ref{030910a}
the S-ring over the above normal $\A$-section is dense. Therefore the S-ring $\A_S$ is dense. Since
it is also primitive by the theorem hypothesis, the number $|S|$ is prime.
Thus in this case statement~(2) holds.\medskip

To complete the proof suppose that $C$ contains no subnormal sections (subnormal flags in
the sense of paper~\cite{EP03be}). Then by Proposition~5.3 of that paper the Cayley scheme~$\CC$
associated with~$\A$ has singularity of degree $\ge 4$ in the pair $(S,T)$ where $S$ and $T$ are the
smallest and greatest sections of the class~$C$. In our terms, this means that $C$ is
a singular class and $|S|=|T|\ge 4$. Then from \cite[Lemma~4.3]{EP03be} it follows that
$\aut(\CC)^S=\sym(S)$ for all $S\in C$.\footnote{This result is a special case of statement~(1)
of Theorem~\ref{260410a} here.} Thus, since $\aut(\A)=\aut(\CC)$, in this case statement~(1) holds.\bull

\section{Generalized wreath product}
\sbsnt{Definition.}\label{100111b}
Let $E_0$ and $E_1$ be equivalence relations on a set~$V$, and $E_0\subset E_1$. Suppose we are given
\nmrt
\tm{a} a group $\Delta_0\le\sym(V/E_0)$ such that the equivalence relation $(E_1)_{V/E_0}$ is
$\Delta_0$-invariant,
\tm{b} for each $X\in V/E_1$ a group $\Delta_X\le\sym(X)$ such that the equivalence relation
$(E_0)_X$ is $\Delta_X$-invariant,
\tm{c} for each $X,X'\in V/E_1$ a bijection $f_{X,X'}:X\to X'$ taking $\Delta_X$ to
$\Delta_{X'}$, and $X/E_0$ to $X'/E_0$, such that
$$
(\Delta_X)^{f_{X,X'}f_{X',X''}}=\Delta_{X''},\qquad X,X',X''\in V/E_1.
$$
\enmrt
The condition in~(c) implies that $\Delta_Xf_{X,X'}=f_{X,X'}\Delta_{X'}=:\Delta_{X,X'}$. Then
$\Delta_X=\Delta_{X,X}$ and $f_{X,X'}\in\Delta_{X,X'}$ for all~$X,X'$.
Thus the data from~(b) and~(c) can be restored from the sets $\Delta_{X,X'}$.
Moreover, instead of the families $\{\Delta_X\}$ from~(b) and $\{f_{X,X'}\}$ from~(c) we could
equivalently be given
\nmrt\label{081110a}
\tm{b-c} for each $X,X'\in V/E_1$ a nonempty set $\Delta_{X,X'}\subset\bij(X,X')$ taking
$(E_0)_X$ to $(E_0)_{X'}$, such that
$$
\Delta_{X,X'}\Delta_{X',X''}=\Delta_{X,X''},\qquad X,X',X''\in V/E_1.
$$
\enmrt
Indeed, set $\Delta_X=\Delta_{X,X}$. Then from the above equality it follows
that $\Delta_X\Delta_X=\Delta_X$, whence by the finiteness of $V$ we conclude
that $\Delta_X$ is a group. Moreover, by the same equality we have
$\Delta_Xf=\Delta_{X,X'}=f\Delta_{X'}$ for any $f\in\Delta_{X,X'}$.\medskip

Set
\qtnl{220110b}
\Gamma=
\{\gamma\in\aut(E_0,E_1):\ \gamma^{V/E_0}\in\Delta_0\ \text{and}\
\gamma^X\in\Delta_{X,X^\gamma},\ X\in V/E_1\}
\eqtn
where $\aut(E_0,E_1)$ is the subgroup of $\sym(V)$ leaving the relations $E_0$ and $E_1$ fixed.
It is easily seen that $\Gamma$ is a group; when $E_0=E_1$, it is permutationally isomorphic to
the wreath product $\Delta_X\wr\Delta_0$ in imprimitive action for any $X\in V/E_1$.
It is also clear that
\qtnl{140111a}
\Gamma_{E_0}=\prod_{X\in V/E_1}\Delta_{X,E_0}
\eqtn
where $\Delta_{X,E_0}=(\Delta_X)_{(E_0)_X}$.

\lmml{230110r}
The equivalence relations $E_0$ and $E_1$ are $\Gamma$-invariant. Moreover, if
\qtnl{220110a}
(\Delta_{X,X'})^{X/E_0,X'/E_0}=(\Delta_0)^{X/E_0,X'/E_0},\quad X,X'\in V/E_1,
\eqtn
then $\Gamma^{V/E_0}=\Delta_0$, and $\Gamma^{X,X'}=\Delta_{X,X'}$ for all $X,X'$.
In particular, the group $\Gamma^{V/E_1}$ is transitive.
\elmm
\proof The first statement as well as the inclusions $\Gamma^{V/E_0}\le\Delta_0$ and
$\Gamma^{X,X'}\subset\Delta_{X,X'}$ follow from the definition of the group~$\Gamma$.
To prove the converse inclusions let $X_0\in V/E_1$. We claim:
{\it given $\delta_0\in\Delta_0$ and $\delta_{X_0}\in\Delta_{X^{}_0,X'_0}$ with
$X'_0$ the class of $E_1$ for which $X'/E_0=(X/E_0)^{\delta_0}$, such that
\qtnl{220710a}
(\delta_{X_0})^{X_0/E_0}=(\delta_0)^{X_0/E_0},
\eqtn
there exists a permutation $\delta\in\Gamma$ such that $\delta^{X_0}=\delta_{X_0}$
and $\delta^{V/E_0}=\delta_0$}. Then $\Gamma^{V/E_0}\ge\Delta_0$
because by ~\eqref{220110a} for any $\delta_0\in\Delta_0$
there exists $\delta_{X_0}\in\Delta_{X^{}_0,X'_0}$ satisfying~\eqref{220710a},
whereas $\Gamma^{X,X'}\supset\Delta_{X,X'}$ follows from the claim
for $X_0=X$ because due to~\eqref{220110a} for any $\delta_X\in\Delta_{X,X'}$
there exists $\delta_0\in\Delta_0$ satisfying~\eqref{220710a}.\medskip

To prove the claim we observe that by~\eqref{220110a}
for each $X\in V/E_1$ other than $X_0$ there exists $\delta_X\in\Delta_{X,X'}$ such that
\qtnl{200710a}
(\delta_X)^{X/E_0}=(\delta_0)^{X/E_0}
\eqtn
where $X'$ is the class of $E_1$ such that $X'/E_0=(X/E_0)^{\delta_0}$. Denote
by $\delta$ the permutation of~$V$ such that $\delta^X=\delta_X$ for
all~$X\in V/E_1$. Then due to~\eqref{220710a} and~\eqref{200710a} we have also
$\delta^{V/E_0}=\delta_0$. Thus $\delta\in\Gamma$.\bull

\dfntnl{260710a}
When (\ref{220110a}) is satisfied, the group $\Gamma$ defined by (\ref{220110b})
is called the generalized wreath product of the family $\{\Delta_{X,X'}\}$ by the
group $\Delta_0$; it will be denoted by $\{\Delta_{X,X'}\}\wr\Delta_0$.
\edfntn

From the arguments in the beginning of the section it follows that the family $\{\Delta_{X,X'}\}$ is uniquely determined by the
group $\Delta_X$ for a fixed $X\in V/E_1$ and the family of bijections $f_{X,X'}:X\to X'$,
$X'\in V/E_1$. Therefore we also can say that the group $\Gamma$ is the generalized
wreath product of the groups $\Delta_X$ and $\Delta_0$ with respect to the above
family of bijections.

\sbsnt{Canonical generalized wreath product.}\label{130111a}
An important special case arises when $V=G$ is a group, the classes of~$E_0$ and $E_1$ are
the left cosets modulo a normal subgroup $L$ of~$G$ and modulo a subgroup $U$ of~$G$ that
contains~$L$ respectively. Suppose that we are also given groups $\Delta_0\le\sym(G/L)$ and
$\Delta_U\le\sym(U)$ such that
\qtnl{100111a}
(G/L)_{right}\le \Delta_0,\quad U_{right}\le \Delta_U,\quad (\Delta_U)^S=(\Delta_0)^S
\eqtn
where $S=U/L$. Given $X,X'\in V/E_1$ set
\qtnl{270910a}
\Delta_{X,X'}=(G_{right})^{X,U}\,\Delta_U\, (G_{right})^{U,X'}.
\eqtn
Then $\Delta_{X,X'}\supset(G_{right})^{X,X'}$ and all conditions in (a) and (b-c) are satisfied: the equivalence relation $(E_1)_{V/E_0}$ is
$\Delta_0$-invariant, $\Delta_{X,X'}$ takes $(E_0)_X$ to $(E_0)_{X'}$ and the equality from~(b-c)
holds. Moreover, from \eqref{100111a} it follows that
$$
(\Delta_{X,X'})^{X/E_0,X'/E_0}=
$$
$$
(G_{right})^{X/E_0,U/E_0}\,(\Delta_U)^{U/L}\, (G_{right})^{U/E_0,X'/E_0}=
$$
$$
((G/L)_{right})^{X/E_0,U/E_0}\,(\Delta_0)^{U/L}\, ((G/L)_{right})^{U/E_0,X'/E_0}=
$$
$$
(\Delta_0)^{X/E_0,X'/E_0}.
$$
Therefore condition~(\ref{220110a}) of Lemma~\ref{230110r} is satisfied and we can form the
generalized wreath product $\Gamma=\{\Delta_{X,X'}\}\wr\Delta_0$.

\dfntnl{110111a}
The group $\Gamma$ is called the canonical generalized wreath product of the group~$\Delta_U$
by the group~$\Delta_0$ (over~$G$); it is denoted by $\Delta_U\wr_S\Delta_0$.
\edfntn

It is easily seen that any canonical generalized wreath product over~$G$
contains~$G_{right}$.

\sbsnt{Automorphism groups.}
Let $\CC=(V,\R)$ be a homogeneous scheme and $E_0,E_1\in\E(\CC)$ such that
$E_0\subset E_1$. Suppose that  $R_{X,Y}=X\times Y$ for all $R\in\R$ contained in $V^2\setminus E_1$
and all $X,Y\in V/E_0$ with $R_{X,Y}\ne\emptyset$. Then we say that the scheme $\CC$
satisfies the {\it $E_1/E_0$-condition} (this definition is obviously equivalent to the
definition given in~\cite{EP03be}).

\thrml{220710c}
Let $\CC$ be a homogeneous scheme on $V$ and $f\in\sym(V)$. Suppose that the $\CC$
satisfies the $E_1/E_0$-condition. Then
$f\in\aut(\CC)$ if and only if
$$
f^{V/E_0}\in\aut(\CC_{V/E_0})\quad\text{and}\quad
f^X\in\iso(\CC_X,\CC_{X'},\varphi_{X,X'}),\ X\in V/E_1,
$$
where $X'=X^f$ and $\varphi_{X,X'}$ is the similarity defined in~\eqref{300810l}.
\ethrm
\proof The family $\{f^X\}$ together with the permutation $f^{V/E_0}$ form
an admissible $E_1/E_0$-pair which is compatible with~$\CC$ in the sense of~\cite{EP03be}.
Thus the required statement is a consequence of Theorem~2.5 of that paper.\bull

\crllrl{050710a}
Let $\CC$ be a homogeneous scheme on $V$ satisfying the $E_1/E_0$-condition. Suppose that the group $\Gamma=\aut(\CC)$ is transitive on $V/E_1$
(this is always true when $\CC$ is a Cayley scheme). Then $\Gamma=\{\Gamma^{X,X'}\}\wr\Gamma^{V/E_0}$.
\ecrllr
\proof Set $\Delta_0=\Gamma^{V/E_0}$ and $\Delta_{X,X'}=\Gamma^{X,X'}$ for all
$X,X'\in V/E_1$. Then obviously the condition from~(a) as well as the condition~\eqref{220110a}
from Lemma~\ref{230110r} are satisfied. Moreover, from the transitivity of $\Gamma$ on $V/E_1$ it follows that the
condition from~(b-c) is also satisfied. Thus we can form the generalized
wreath product $\Gamma'=\{\Delta_{X,X'}\}\wr\Delta_0$. Clearly, $\Gamma'\ge\Gamma$.
The converse inclusion follows from Theorem~\ref{220710c}.\bull\medskip

Let us consider an example of the situation from Corollary~\ref{050710a}.
Let $\CC=(V,\R)$ be a homogeneous scheme satisfying the $E_1/E_0$-condition. Set
$$
\Delta_0=\aut(\CC_{V/E_0}),\qquad\Delta_{X,X'}=\iso(\CC_{X^{}},\CC_{X'},\varphi_{X^{},X'}),\quad
X,X'\in V/E_1.
$$
In general, we cannot form the generalized wreath product $\{\Delta_{X,X'}\}\wr\Delta_0$
because the condition from (b-c) or equality~\eqref{220110a} is not necessarily true.
Now suppose  that {\it the scheme $\CC_{V/E_0}$ and all schemes $\CC_X$ with $X\in V/E_1$
are regular} (i.e. schemes of regular permutation groups). Then both of these conditions follow
from the fact that any similarity from regular scheme to another scheme is induced by an
isomorphism. Thus in this case the above generalized wreath product can be constructed. But
then
$$
\aut(\CC)=\{\Delta_{X,X'}\}\wr\Delta_0
$$
by Theorem~\ref{220710c}.

\sbsnt{Schurity and non-schurity of S-rings.}\label{241210a}
In the rest of the section we are going to get a necessary and sufficient condition
for an S-ring over an abelian group that is a generalized wreath product, to be
schurian. This condition as well as a criterion for non-schurity will be deduced
from the following result.

\thrml{260710c}
Let $\CC$ be a homogeneous scheme on~$V$ satisfying the $E_1/E_0$-condition, and
$\Gamma=\{\Delta_{X,X'}\}\wr\Delta_0$ a generalized wreath product 
such that $\Delta_0\le\aut(\CC_{V/E_0})$ and $\Delta_{X,X'}\subset\iso(\CC_X,\CC_{X'},\varphi_{X,X'})$
for all $X,X'\in V/E_1$. Suppose that
\qtnl{260710d}
\orb(\Delta_{X,E_0})=X/E_0,\qquad X\in V/E_1.
\eqtn
Then $\CC=\inv(\Gamma)$
if and only if $\CC_{V/E_0}=\inv(\Delta_0)$ and $\CC_X=\inv(\Delta_X)$ for all $X\in V/E_1$.
\ethrm
\proof The necessity follows from Lemma~\ref{230110r} because
$\inv(\Gamma)^{V/E_0}=\inv(\Gamma^{V/E_0})$ and $\inv(\Gamma)^X=\inv(\Gamma^X)$
for all $X\in V/E_1$. To prove the sufficiency we check
that each relation $R\in\R(\CC)$ is an orbit of the group $\Gamma$. We note that
since $\Gamma^{V/E_0}=\Delta_0$ (Lemma~\ref{230110r}) and $\CC_{V/E_0}=\inv(\Delta_0)$, the group $\Gamma^{V/E_0}$
acts transitively on $R_{V/E_0}$. On the other hand, since $\Delta_{X,X'}\subset\iso(\CC_X,\CC_{X'},\varphi_{X,X'})$
and the scheme~$\CC$ satisfies the $E_1/E_0$-condition, the
group~$\Gamma$ acts on the non-empty sets $R_{X,X'}$, $X,X'\in V/E_1$. Thus $\Gamma$ acts on~$R$. To prove that
this action is transitive, first suppose that $R\subset E_1$. Then $R$ is a disjoint union of
$R_X$, $X\in V/E_1$. Since the group $\Gamma^{V/E_1}$ is transitive (Lemma~\ref{230110r}), it
suffices to verify that the group $\Gamma^X$ acts transitively on each~$R_X$. However, this is true
because $\Gamma^X=\Delta_X$ (Lemma~\ref{230110r}) and $\CC_X=\inv(\Delta_X)$.\medskip

Now let $R\subset V^2\setminus E_1$.  Then given $Y,Y'\in V/E_0$ we have either $R_{Y^{},Y'}=\emptyset$
or $R_{Y^{},Y'}=Y\times Y'$. So it suffices to verify that in the latter case the group
$\Gamma_{Y^{},Y'}=\Gamma_{\{Y^{}\}}\cap \Gamma_{\{Y'\}}$ acts transitively on $R_{Y^{},Y'}$.
However,
by~\eqref{140111a}  we have
$$
\Delta_{X^{},E_0}\times\Delta_{X',E_0}\times\{\id_{V\setminus(X\cup X')}\}\le\Gamma_{E_0}\le\Gamma_{Y^{},Y'}
$$
where $X$ and $X'$ are the classes of the equivalence relation $E_1$ containing $Y$ and $Y'$
respectively. Thus the required statement follows from~\eqref{260710d}.\bull

\crllrl{260210a}
Let $\A$ be an S-ring over an abelian group~$G$. Suppose that $\A=\A_U\wr_S\A_{G/L}$
for some $\A$-section~$S=U/L$.
Then $\A$ is schurian if and only if so are the S-rings $\A_{G/L}$ and $\A_U$
and there exist groups $\Delta_0\le\sym(G/L)$ and $\Delta_U\le\sym(U)$
satisfying~\eqref{100111a} and such that
\qtnl{100111c}
\Delta_0\twoe\aut(\A_{G/L})\quad\text{and}\quad\Delta_U\twoe\aut(\A_U).
\eqtn
Moreover, in this case $\aut(\A)\twoe\Delta_U\wr_S\Delta_0$.
\ecrllr
\proof It is easily seen that the S-ring $\A$ satisfies the $U/L$-condition if and only if
the Cayley scheme $\CC$ associated with~$\A$ satisfies the $E_1/E_0$-condition where
$E_0=E_L$ and $E_1=E_U$. Thus the necessity follows from Corollary~\ref{050710a}
with $\Delta_0=\aut(\A)^{G/L}$ and $\Delta_U=\aut(\A)^U$.\medskip

To prove the sufficiency suppose we are given groups $\Delta_0\le\sym(G/L)$ and
$\Delta_U\le\sym(U)$ satisfying~\eqref{100111a} and~\eqref{100111c}. Set
$\Gamma=\Delta_U\wr_S\Delta_0$ (as in the Subsection~\ref{130111a})
and verify that the hypothesis of Theorem~\ref{260710c} is satisfied.
Indeed, the inclusion $\Delta_0\le\aut(\CC_{G/L})$ is clear
whereas the inclusion $\Delta_{X,X'}\subset\iso(\CC_X,\CC_{X'},\varphi_{X,X'})$ where
$X,X'\in V/E_1$, is true due to~\eqref{270910a} because
$(G_{right})^{Y,Y'}\subset\iso(\CC_Y,\CC_{Y'},\varphi_{Y,Y'})$ for all $Y,Y'\in V/E_1$.
Moreover, from the definition of $\Delta_{X,X'}$ it follows that
$\Delta_X=(\Delta_U)^f$ where $f\in (G_{right})^{U,X}$. Therefore
by the second inclusion of~\eqref{100111a} this implies that
$$
\orb(\Delta_{X,E_0})=
\orb((\Delta_{U,E_0})^f)=
\orb(\Delta_{U,E_0})^f=
(U/L)^f=X/E_0.
$$
So condition~\eqref{260710d} of Theorem~\ref{260710c} is also satisfied. Thus
by this theorem we conclude that $\CC=\inv(\Gamma)$ whence it follows that the
S-ring~$\A$ is schurian.\bull\medskip

The following statement gives a criterion for the non-schurity of an S-ring over an
abelian group.

\crllrl{250510a}
Let $\A$ be an S-ring over an abelian group~$G$. Suppose that $\A=\A_U\wr_S\A_{G/L}$
for some $\A$-section~$S=U/L$. Then the S-ring $\A$ is non-schurian whenever
\qtnl{250510c}
\aut(\A_U)^S\cap\aut(\A_{G/L})^S\not\twoe\aut(\A_S).
\eqtn
\ecrllr
\proof Suppose on the contrary that the S-ring $\A$ is schurian. Then by
Corollary~\ref{260210a} the S-rings $\A_{G/L}$ and $\A_U$ are schurian
and there exist groups $\Delta_0\le\sym(G/L)$ and $\Delta_U\le\sym(U)$
satisfying conditions~\eqref{100111a} and~\eqref{100111c}. Therefore
$(\Delta_0)^S\le\aut(\A_S)$, $(\Delta_U)^S\le\aut(\A_S)$ and
$$
(\Delta_0)^S\twoe\aut(\A_{G/L})^S\twoe\aut(\A_S).
$$
Thus the intersection in the left-hand side of (\ref{250510c}) contains the
subgroup $(\Delta_U)^S=(\Delta_0)^S$ which is $2$-equivalent to the
group $\aut(\A_S)$. Contradic\-tion.\bull

\section{Isolated classes}\label{041010b}
\sbsnt{Isolated pairs.}\label{081110d}
Let $S=L_1/L_0$ and $T=U_1/U_0$ be nontrivial sections of an S-ring~$\A$ over an abelian 
group~$G$.

\dfntnl{1111111}
We say that $S$ and $T$ form an {\it isolated pair} in~$\A$ if $T$ is a multiple of $S$ and conditions~(S1)
and~(S2) are satisfied.
\edfntn

It immediately follows from the definition that $U_0$, $U_1\setminus U_0$ and $G\setminus U_1$
are $\A$-subsets of~$G$. Moreover, the set $\S=\S(\A)$ is uniquely determined by the
sets $\S_{U_0}=\S(\A_{U_0})$, $\S_{L_1/L_0}=\S(\A_{L_1/L_0})$ and $\S_{G/L_1}=\S(\A_{G/L_1})$ as follows:
\qtnl{220410b}
\S=\S_{U_0}\,\cup\, (\pi_0^{-1}({\S_{L_1/L_0}}))_{L_1\setminus L_0}\S_{U_0}\,\cup\,
(\pi_1^{-1}(\S_{G/L_1}))_{G\setminus U_1},
\eqtn
where $\pi_0=\pi_{L_1/L_0}$ and $\pi_1=\pi_{G/L_1}$ are the natural epimorphisms.
Moreover, the three sets on the right-hand side are pairwise disjoint and equal to
$\S(\A)_{U_0}$, $\S(\A)_{U_1\setminus U_0}$ and $\S(\A)_{G\setminus U_1}$
respectively.\medskip

Obviously, the $\A$-sections forming an isolated pair are projectively
equivalent. The projective equivalence class containing them will be
called {\it isolated}; we also say that it contains the pair.

\lmml{260410o}
Let $\A$ be a circulant S-ring. Then any isolated class $C\in\P(\A)$
contains exactly one isolated pair. This pair consists of the smallest and the largest elements
of~$C$.
\elmm
\proof By Theorem~\ref{281010a} there exist the smallest and greatest sections
$L_1/L_0$ and $U_1/U_0$ in the class~$C$. Clearly $\pi(L_1)$ and $\pi(U_0)$ are
$\A_{\pi(U_1)}$-subgroups where $\pi=\pi_{U_1/L_0}$. Moreover, since $U_1/U_0$ is a multiple
of $L_1/L_0$ we also have
\qtnl{291010c}
\pi(U_1)=\pi(L_1)\times\pi(U_0).
\eqtn
Let $L'_1/L'_0$ and $U'_1/U'_0$ be sections in~$C$ forming an isolated pair in $\A$.
Then it suffices to verify that $U'_1=U_1$ and $L'_0=L_0$ (then obviously $U'_0=U_0$
and $L'_1=L_1$). To do this we observe that by the definition of the smallest and largest
sections we have
\qtnl{060510a}
L_1\le L'_1\le U'_1\le U_1\quad\text{and}\quad
L_0\le L'_0\le U'_0\le U_0.
\eqtn
Suppose first that $U'_1\ne U_1$. Then there exists $X\in\S(\A)_{U_1\setminus U'_1}$.
From the second equality of~(S1) with $U_1=U'_1$ and $L_1=L'_1$ it follows that
$L'_1X=X$, and hence due to left-hand side of~(\ref{060510a}) also $L_1X=X$. This implies that
\qtnl{291010a}
\pi(L_1)\pi(X)=\pi(X).
\eqtn
On the other hand,  due to~\eqref{291010c} Lemma~\ref{130209d} implies that
$\pr_{\pi(L_1)}(\pi(X))$ is a basic set of the S-ring $\A_{\pi(L_1)}$.
However, by~\eqref{291010a} this set coincides with~$\pi(L_1)$, which
is impossible because $\pi(L_1)=L_1/L_0\ne\{1\}$. Thus $U'_1=U_1$.\medskip

To complete the proof suppose that $L_0\ne L'_0$. Then there exists a set
$X'\in\S(\A)_{L'_0\setminus L_0}$. By~\eqref{060510a} we have $X'\subset U_0$
and hence $\pi(X')\subset\pi(U_0)$. Due to~\eqref{291010c} and the inequality  $L_1\ne L_0$,
the full $\pi$-preimage
of $\pi(X')$ does not coincide with $X'$. Therefore one can find (in this preimage) a basic set
$X\in\S(\A)_{U_1\setminus U_0}$ such that $\pr_{\pi(U_0)}(\pi(X))=\pi(X')$. On the other
hand, from the first equality of~(S1) with $U_0=U'_0$ and $L_0=L'_0$ it follows that
$L'_0X=X$. This implies that $\pi(L'_0)=\pr_{\pi(U_0)}(\pi(X))=\pi(X')$ is a basic
set of the S-ring $\A_{\pi(U_0)}$. However, this is impossible because
$\pi(X')\ne\{1\}$. Thus $L'_0=L_0$.\bull\medskip

Below for an S-ring $\A$ over a cyclic group the smallest and the largest elements
of a class $C\in\P(\A)$ the existence of which follows from Theorem~\ref{281010a}, are
denoted by $S_{min}(C)=L_1(C)/L_0(C)$ and $S_{max}(C)=U_1(C)/U_0(C)$ respectively.

\crllrl{130710a}
Let $C\in\P(\A)$ be a primitive isolated class of a circulant S-ring $\A$.
Suppose that an $\A$-section $S=U/L$ has no subsection from this class.
Then either $L\ge L_1$ or $U\le U_0$ where $L_1=L_1(C)$ and $U_0=U_0(C)$.
\ecrllr
\proof By Lemma~\ref{260410o} the sections $S_{min}=L_1/L_0$ and $S_{max}=U_1/U_0$
where $L_0=L_0(C)$ and $U_1=U_1(C)$, form an isolated pair in~$\A$. Suppose that
the section $S$ is such that $L\not\ge L_1$ and $U\not\le U_0$. Then
\qtnl{130710b}
U\ge L_1\quad\text{and}\quad L\le U_0.
\eqtn
Indeed, since the S-ring $\A$ satisfies the $U_1/L_1$-condition, it follows that either
$U\ge L_1$, or $U\le U_1$. In the latter case the right-hand side of~\eqref{130710b}
is obvious, whereas the left one follows from the assumption $U\not\le U_0$, the primitivity 
of~$C$ and the equality $L_1\cap U_0=L_0$. The former case is proved in a similar way. 
\medskip

Now from the left-hand side of~\eqref{130710b} it follows that
$$
S=U/L\ge U\cap U_1/L L_0\ge L L_1/L L_0.
$$
Besides, the latter section is projectively equivalent to $S_{min}$ because
by the right-hand side of~\eqref{130710b} we have $L_1\cap LL_0=L_0$, and $L_1\,LL_0=L L_1$.
Thus $S$ has a subsection from~$C$. Contradiction.\bull\medskip

\sbsnt{Extension construction.}
Let $\A$ be an S-ring over a cyclic group $G$ and $C\in\P(\A)$ an isolated class.
Then by definition~\ref{1111111} and Lemma~\ref{260410o} the S-ring~$\A$ satisfies 
the $U_i/L_i$-condition where $U_i=U_i(C)$ and $L_i=L_i(C)$,
$i=0,1$, and $\A_{U_1/L_0}=\A_{L_1/L_0}\otimes\A_{U_0/L_0}$.\medskip 

Suppose we are additionally 
given an S-ring $\B$ over the group $S=L_1/L_0$, such that $\B\ge\A_S$.
Then by Theorem~\ref{160710a} there are uniquely determined S-rings
\qtnl{301110b}
\A_1=\A_{U_0}\wr_{U_0/L_0}(\B\otimes\A_{U_0/L_0})\quad\text{and}\quad
\A_2=(\B\otimes\A_{U_0/L_0})\wr_{U_1/L_1}\A_{G/L_1}
\eqtn
over the groups $U_1$ and $G/L_0$ respectively. Obviously, the restrictions of these S-rings
to $U_1/L_0$ coincide with $\B\otimes\A_{U_0/L_0}$. Therefore by Theorem~\ref{160710a} there
is a uniquely determined S-ring
\qtnl{091110a}
\Ext_C(\A,\B)=\A_1\wr_{U_1/L_0}\A_2
\eqtn
over the group~$G$. 

\dfntnl{2222222}
The S-ring $\Ext_C(\A,\B)$ is called the extension of the S-ring~$\A$ by means of
the S-ring~$\B$ with isolated class~$C$.
\edfntn

The following statement is straightforward.

\lmml{011110a}
In the above notation set $\A'=\Ext_C(\A,\B)$. Then
\nmrt
\tm{1} $\A'\ge\A$; moreover, $\A=\A'$ if and only if $\B=\A_S$,
\tm{2} $\A'$ satisfies the $U_0/L_0$- and $U_1/L_1$-conditions,
\tm{3} $\A'_{U_0}=\A^{}_{U_0}$, $\A'_{G\setminus L_1}=\A^{}_{G\setminus L_1}$ and
$\A'_{U_1/L_0}=\B\otimes\A^{}_{U_0/L_0}$.\bull
\enmrt
\elmm

From statements~(2) and~(3) of Lemma~\ref{011110a} it follows that the $\A'$-sections
$S_{min}(C)$ and $S_{max}(C)$
form an isolated pair in the S-ring~$\A'$. Therefore by Lemma~\ref{260410o} we
obtain the following statement.

\thrml{160410a}
Let $C\in\P(\A)$ be an isolated class of a circulant S-ring $\A$ and
$\A'=\Ext_C(\A,\B)$ for some $\B$ as above. Then the class $C'\in\P(\A')$ containing $S=S_{min}(C)$
and $T=S_{max}(C)$, is isolated and $S_{min}(C')=S$, $S_{max}(C')=T$.\bull
\ethrm

The following statement gives a necessary and sufficient condition for the schurity of
extension~\eqref{091110a}.

\thrml{301110a}
Let $C\in\P(\A)$ be an isolated class of a schurian circulant S-ring $\A$.
Suppose that the S-rings $\A_S$ and $\B$ are schurian. Then the S-ring
$\A'=\Ext_C(\A,\B)$ is schurian if and only if the S-ring $\A$ is schurian.
\ethrm
\proof Suppose first that the S-ring $\A$ is schurian. Then the S-ring $\A_{U_0}$
is schurian, and due to the schurity of~$\B$ also the S-ring $\B\otimes\A_{U_0/L_0}$.
Therefore by Corollary~\ref{260210a} applied to the S-ring~$\A_1$
defined in \eqref{301110b} and the groups $U=U_0$, $\Delta_0=\aut(\B)\times\aut(\A)^{U_0/L_0}$,
$\Delta_U=\aut(\A)^{U_0}$, this S-ring is schurian and
\qtnl{301110c}
\aut(\A_1)^{U_1/L_0}=\aut(\B)\times\aut(\A)^{U_0/L_0}.
\eqtn
Similarly, the S-ring~$\A_2$ defined in \eqref{301110b} is schurian and
\qtnl{301110d}
\aut(\A_2)^{U_1/L_0}=\aut(\B)\times\aut(\A)^{U_0/L_0}.
\eqtn
Thus, the S-ring $\A'$ is the generalized wreath product of two schurian S-rings $\A_1$ and $\A_2$
(see~\eqref{091110a}). Thus by \eqref{301110c} and \eqref{301110d} Corollary~\ref{260210a}
with $\Delta_0=\aut(\A_2)$ and $\Delta_U=\aut(\A_1)$ implies that $\A'$ is schurian.
The converse statement is proved analogously with taking into account equalities
from statement~(3) of Lemma~\ref{011110a}.\bull

\sbsnt{Automorphisms.}\label{161110x}
Let $S=L_1/L_0$ and $T=U_1/U_0$ be sections of the group $G$ such that $T$ is a multiple of~$S$.
For a transitive group $M\le\sym(S)$ and cosets $X,X'\in T$ set
\qtnl{231110a}
\Delta_0=M\times\{\id_{U_0/L_0}\}\quad\text{and}\quad
\Delta_{X,X'}=(((U_1)_{right})_{E_{L_1}})^{X,X'}.
\eqtn
Then obviously the equivalence relation $(E_1)_{U_1/E_0}$ where $E_1=E_{U_0}$ and $E_0=E_{L_0}$,
is $\Delta_0$-invariant. Moreover, all conditions in (a) and (b-c) on page~\pageref{081110a} are
satisfied for $V=U_1$. Finally, given $X,X'\in U_1/E_1$ there exist $a,a'\in L_1$ such that
$X=aU_0$ and $X'=a'U_0$ (we note that $a$ and $a'$ are uniquely determined modulo $L_0$).
Then from the transitivity of~$M$ and the definition of $\Delta_0$ it follows that
$$
(\Delta_{X,X'})^{X/E_0,X'/E_0}=\{\gamma_{a,a'}\}=(\Delta_0)^{X/E_0,X'/E_0}
$$
where $\gamma_{a,a'}$ is the bijection from $X/E_0$ onto $X'/E_0$ taking
$axL_0$ to $a'xL_0$, $x\in U_0$. Thus condition~\eqref{220110a} of Lemma~\ref{230110r} is
also satisfied and one can consider the generalized wreath product
$\{\Delta_{X,X'}\}\wr\Delta_0$.\medskip

Denote by $\gwr(S,T,M)$ the subgroup $\Delta$ of the pointwise
stabilizer of the set $G\setminus U_1$ in the group $\sym(G)$ such that
\qtnl{260410c}
\Delta^{U_1}=\{\Delta_{X,X'}\}\wr\Delta_0.
\eqtn
\noindent It is easily seen that given a group $H\le G$ such
that $(L_1\cap L_0H)/L_0$ is a block of~$M$ we have
$$
\gwr(S,T,M)^{\pi(G)}=\gwr(\pi(S),\pi(T),M^{\pi(S)})
$$
where $\pi=\pi_{G/H}$.\medskip

When $\A$ is a circulant S-ring, $C\in\P(\A)$ and
$S=S_{min}(C)$, $T=S_{max}(C)$, we set $\gwr_\A(C,M)=\gwr(S,T,M)$. The first statement of
the following theorem generalizes~\cite[Lemma~4.3]{EP03be}.

\thrml{260410a}
Let $\A$ be a circulant S-ring, $C\in\P(\A)$ an isolated class,
and $\Delta=\gwr_\A(C,M)$ where $M=\aut(\A_S)$ with $S=S_{min}(C)$.
Then
\nmrt
\tm{1} $\Delta\le\aut(\A)$ ,
\tm{2}  $\Delta^{S'}=\aut(\A_{S'})$ for all $S'\in C$,
\tm{3} $\Delta^{S'}\le (S')_{right}$ whenever $C$ is primitive and
a section $S'\in\scc(\A)\setminus C$ is either primitive, or of order coprime to~$|S|$.
\enmrt
\ethrm
\proof Let $f\in\Delta$. Then $f^{G/U_1}=\id_{G/U_1}$ and $f^X=\id_X$ for all cosets
$X\in G/U_1$ other than $U_1$ (we keep the notations as
in Subsection~\ref{081110d}). Since $\A$ satisfies the $U_1/L_1$-condition, by
Theorem~\ref{220710c} with $\CC$ being the Cayley scheme associated with~$\A$, it suffices
to verify that $g:=f^{U_1}$ is an automorphism of the S-ring~$\A_{U_1}$. However,
from the definition of the group~$\Delta$ it follows that
$$
g^X\in ((U_1)_{right})^{X,X'}\subset\iso(\CC_X,\CC_{X'},\varphi_{X,X'})
$$
for all $X\in U_1/U_0$ where $X'=X^g$, and
$$
g^{U_1/L_0}\in \Delta_0=M\times\id_{U_0/L_0}\le\aut(\CC_{U_1/L_0})
$$
because $\A_{U_1/L_0}=\A_{L_1/L_0}\otimes\A_{U_0/L_0}$.
Since the S-ring $\A_{U_1}$ satisfies the $U_0/L_0$-condition, we obtain
by Theorem~\ref{220710c} that $g\in\aut(\A_{U_1})$.\medskip

To prove statement~(2) let $S'\in C$. Denote by $\Delta'$ the subgroup of $\sym(U_1)$
generated by $\Delta^{U_1}$ and $(U_1)_{right}$. By statement~(1) we have
\qtnl{240111a}
\Delta'\le\aut(\A)^{U_1}.
\eqtn
So $(\Delta')^{S'}\le\aut(\A)^{S'}$. Moreover, by Lemma~\ref{160210a} the groups
$(\Delta')^{S'}$ and $(\Delta')^{S^{}}$ are isomorphic. On the other hand,
since $(U_1)_{right}$ normalizes $\Delta^{U_1}$ and $\Delta^S\ge S_{right}$, we have
$(\Delta')^S=\Delta^S=\aut(\A_S)$. Thus $(\Delta')^{S'}=\aut(\A_{S'})$ by 
Theorem~\ref{261010a}.\medskip

To prove statement~(3) let $S'=U'/L'$ is either a primitive section not in~$C$, or a section 
of order coprime to~$|S|$. Then obviously this section has no subsection from $C$.
So by Corollary~\ref{130710a} with $S=S'$ we have either $L'\ge L_1$ or $U'\le U_0$. To 
complete the proof it suffices to note that in the former case $\Delta^{S'}=\id_{S'}$ 
whereas in the latter case $\Delta^{S'}=(S')_{right}$.\bull

\section{Extension of a singular S-ring}\label{041010c}

Any singular class defined in Subsection~\ref{170111a} is obviously isolated. The following result
shows how the set $P_{sgl}(\A)$ of all singular classes in the S-ring~$\A$ varies when passing
to a special case of extension~\eqref{091110a}.

\thrml{200410f}
Let $\A$ be a circulant S-ring, $C\in\P_{sgl}(\A)$ a singular class of prime
order and $\A'=\Ext_C(\A,\Z S)$ where $S=S_{min}(C)$. Then
\nmrt
\tm{1} $\H(\A)=\H(\A')$; in particular, $\scc(\A)=\scc(\A')$ and $\P(\A)=\P(\A')$,
\tm{2} $\P_{sgl}(\A')=\P_{sgl}(\A)\setminus\{C\}$.
\enmrt
\ethrm
\proof We keep the notations as in Subsection~\ref{081110d}.
To prove statement (1) it suffices to verify that any $\A'$-subgroup~$H$ belongs
to~$\H(\A)$. However, if $H\not\subset U_1$, then $H$ is generated by an element in $G\setminus U_1$.
Denote by $X$ the basic set of $\A'_H$ containing this element. Then $H=\lg X\rg$. On the
other hand, $X\in\S(\A)$ by statements~(2) and~(3) of Lemma~\ref{011110a}. Thus $H\in\H(\A)$.
Let $H\subset U_1$. By statement~(3) of the same lemma
we have $\H(\A_{U_0})=\H(\A'_{U_0})$. So we can assume that $H\not\subset U_0$. Then $H$ contains
a set from $\S(\A')_{U_1\setminus U_0}$. Since by statement~(2) of the same lemma the
S-ring $\A'$ satisfies the $U_0/L_0$-condition, the group generated by this set, and hence
the group $H$, contains~$L_0$. Taking into account that $H\not\subset U_0$ and the fact
that the group $\pi(U_1)$  with $\pi=\pi_{G/L_0}$ is the direct product of the group $L_1/L_0$
of prime order and the group $\pi(U_0)$, we conclude that $H=L_1H'$ where $H'$ is the full
$\pi$-preimage of the group $\pi(H)\cap \pi(U_0)$. Since $L_1$ and $H'$ are $\A$-subgroups,
we are done.\medskip

To prove statement (2) we observe that since the order of $C$ is at least~$3$, then
$\rk(\A'_S)>2$. Therefore the class $C$ is not singular in $\A'$. For the rest of the
proof we need the following auxiliary lemma. Below we set $p=|S|$.

\lmml{091110b}
Let $X\in\S(\A)$ and $X'\in\S(\A')$ be such that $X'\subset X$. Then
$\rad(X)=\rad(X')$. Moreover, if $X\subset U_1\setminus U_0$, then
\qtnl{110510c}
X=\bigcup_{\sigma\in T_p} (X')^\sigma\quad\text{and}\quad|X|=(p-1)|X'|
\eqtn
where $T_p$ is the subgroup of $\aut(G)$ of order $p-1$ that acts trivially on each
Sylow $q$-subgroup of $G$, $q\ne p$.
\elmm
\proof If $X\subset U_0$ or $X\subset G\setminus U_1$, then the required statement immediately
follows from statement~(3) of Lem\-ma~\ref{011110a}. Suppose that $X\subset U_1\setminus U_0$. Then by (\ref{220410b})
there exists $Y\in\S(\A)_{U_0}$ and $x'\in L_1\setminus L_0$ such that
\qtnl{110510a}
X=(L_1\setminus L_0)Y,\qquad X'=(x'L_0)Y.
\eqtn
Since $Y\subset U_0$ and $U_0\cap L_1=L_0$, we have
$$
\rad(X)=\rad(L_1\setminus L_0)\rad(Y)=L_0\rad(Y).
$$
On the other hand, obviously $\rad(X')=L_0\rad(Y)$. Thus, $\rad(X)=\rad(X')$.
Since $\pi(U_1)=\pi(L_1)\times \pi(U_0)$ where $\pi=\pi_{G/L_0}$, the group $T_p$
leaves the set $L_0Y$ fixed. Therefore (\ref{110510c}) follows from~(\ref{110510a}).\bull\medskip

Let us prove that any class $\wt C\ne C$ belonging to the set $\P(\A)=\P(\A')$
is singular in $\A$ and $\A'$ simultaneously. Set
$$
\wt L_i=L_i(\wt C)\quad\text{and}\quad\wt U_i=U_i(\wt C),\qquad i=0,1
$$
(see Subsection~\ref{081110d}). Then by Theorem~\ref{260410o}
it suffices to verify that conditions~(S1) and~(S2)
for $L_i=\wt L_i$ and $U_i=\wt U_i$ are satisfied for $\A$ and $\A'$ simultaneously, and
that $\rk(\A_{H_1})=2$ if and only if $\rk(\A'_{H_1})=2$ where $H_1=\wt L_1/\wt L_0$ and
$H_2=\wt U_0/\wt L_0$. However, from the first statement of Lemma~\ref{091110b} it follows that
the S-rings $\A$ and $\A'$ satisfy the $\wt U_0/\wt L_0$-condition (or the $\wt U_1/\wt L_1$-condition)
simultaneously. Thus it suffices to verify that
$$
\A_H=\A_{H_1}\otimes \A_{H_2},\ \rk(\A_{H_1})=2
\quad\Leftrightarrow\quad
\A'_H=\A'_{H_1}\otimes \A'_{H_2},\ \rk(\A'_{H_1})=2
$$
where $H=\wt U_1/\wt L_0$.  To do this we note that $H=H_1\times H_2$ and
$H_i$ is a $\A^{}_H$- and $\A'_H$-subgroups. By Lemma~\ref{130209d} this implies
that given $X\in\S(\A)$ (resp. $X'\in\S(\A')$) the set $\wt X_i=\pr_{H_i}\wt\pi(X)$
(resp. $\wt X'_i=\pr_{H_i}\wt\pi(X')$) is an $\A$-group (resp. $\A'$-group) for
$i=0,1$, where $\wt\pi=\pi_{\wt U_1/\wt L_0}$.
Thus the required statement is a consequence of the following lemma.

\lmm
Suppose that $\rk(\A_{H_1})=2$. Then given $X\in\S(\A)$ and $X'\in\S(\A')$ such that
$X'\subset X\subset\wt U_1$, we have
\qtnl{130510a}
\wt X=\wt X_1\times\wt X_2
\quad\Leftrightarrow\quad
\wt X'=\wt X'_1\times\wt X'_2
\eqtn
where $\wt X=\wt\pi(X)$ and $\wt X'=\wt\pi(X')$. Moreover, in any case $\wt X^{}_1=\wt X'_1$.
\elmm
\proof Without loss of generality we can assume that $\wt X\ne\wt X'$.
Then obviously $X\ne X'$. Therefore $X\subset U_1\setminus U_0$ by
statement~(3) of Lemma~\ref{011110a}. Due to~\eqref{110510a} this implies that
$L_1\subset\lg X\rg\subset\wt U_1$. We claim that
\qtnl{130510u}
\wt\pi(L_1)\subset H_2.
\eqtn
Suppose on the contrary that this is not true. Denote by $M_i$ (resp. $N_i$) the $H_1$-projection
(resp. the $H_2$-projection) of the group $\wt\pi(L_i)$, $i=0,1$. Then
$\wt\pi(L_i)=M_i\times N_i$. Moreover, taking into account that
$\wt\pi(L_i)$ is an $\A$-group, we conclude by Lemma~\ref{130209d} that
$M_i$ and $N_i$ are $\A$-groups. Since also $\rk(\A_{H_1})=2$, the group $M_i$ is either $1$ or
$H_1$, and either $M_1=M_0$ or $N_1=N_0$. In addition, our supposition implies that $M_1=H_1$. 
\medskip

Next, the S-ring $\A$ satisfies the $U_0/L_0$-condition and
$\rk(\A_{L_1/L_0})=2$. Therefore $L_1\setminus L_0\in\S(\A)$, and hence the set
$\wt Z=\wt\pi(L_1\setminus L_0)$ is a basic one. Besides,
\qtnl{151110a}
\wt Z=\css
\wt\pi(L_1),                        &\text{if $\wt\pi(L_1)=\wt\pi(L_0)$,}\\
\wt\pi(L_1)\setminus\wt\pi(L_0),    &\text{otherwise.}\\
\ecss
\eqtn
(Indeed, in the first case this is true because the group $\wt\pi(L_1)$ is the
union of the sets $\wt\pi(xL_0)=\wt\pi(x)\wt\pi(L_0)$, $x\in L_1$; in the second case the basic
set $\wt Z$ obviously contains the $\A$-set $\wt\pi(L_1)\setminus\wt\pi(L_0)$ and hence
coincides with it.) Now, the equality $\wt Z=\wt\pi(L_1)$ is impossible because
otherwise the set $H_1=\pr_{H_1}(\wt Z)$ is a basic one which is possible only if $H_1=1$.
Thus by~\eqref{151110a} we have
$$
\wt Z=\wt\pi(L_1)\setminus\wt\pi(L_0)
$$
and so either $\pr_{H_1}(\wt Z)=H_1$ or $\pr_{H_2}(\wt Z)=N_1$. As before the former case
is impossible because $\wt Z$ is a basic set, and the latter case is possible only if $N_1=1$.
Thus
$$
\wt\pi(L_0)=1\quad\text{and}\quad\wt\pi(L_1)=H_1.
$$
It follows that $\wt L_1=L_1\wt L_0$ and $\wt L_0\ge L_0$. Therefore $L_0\le L_1\cap \wt L_0<L_1$.
Since $\rk(\A_{L_1/L_0})=2$ is prime, this implies that $L_0=L_1\cap \wt L_0$. Thus
the section $\wt L_1/\wt L_0\not\in C$ is a multiple of the section $L_1/L_0\in C$.
Contradiction. Thus claim~\eqref{130510u} is proved.\medskip

Let us complete the proof of (\ref{130510a}). First, since $X\subset U_1\setminus U_0$
(see above), from~(\ref{110510a}) it follows that $\wt X=\wt\pi(L_1\setminus L_0)\wt Y$ and
$\wt X'=\wt\pi(x'L_0)\wt Y$ where $\wt Y=\wt\pi(Y)$. Due to~(\ref{130510u}) this implies that
\qtnl{130510g}
\wt X^{}_1=\wt X'_1=\wt Y_1
\eqtn
where $\wt Y_1=\pr_{H_1}\wt Y$. Next, taking into account that the radicals of the
sets $\wt X$ and $\wt X'$ contains the group $\wt\pi(L_0)$, without loss of generality
we can assume that this group is trivial. Therefore by Lemma~\ref{091110b} we have
\qtnl{130510h}
\wt X_2=\bigcup_{\sigma\in T_p}\wt\pi(x')^\sigma\wt Y_2,\quad
\wt X'_2=\wt\pi(x')\wt Y_2.
\eqtn
If $\wt\pi(L_1)=1$, then $\wt X_2=\wt X'_2=\wt Y_2$ and the required statement immediately
follows from~\eqref{130510g}. Suppose that the group $\wt\pi(L_1)$ is nontrivial. Then it is an $\A$-group
of order~$p$ and the set $\wt\pi(L_1)\setminus\{1\}$ is an orbit of the group~$T_p$ that contains
the element $\wt\pi(x')$. Moreover, by condition~(S2) we have $\wt\pi(L_1)\lg\wt Y\rg=\wt\pi(L_1)\times\lg\wt Y\rg$.
Therefore the union in~\eqref{130510h} is a disjoint one. Thus to complete the proof of
equivalence~\eqref{130510a} it suffices to see that due to equalities~\eqref{130510g}
and~\eqref{130510h} both $\wt X=\wt X_1\times\wt X_2$ and $\wt X'=\wt X'_1\times\wt X'_2$
are equivalent to the equality $\wt Y=\wt Y_1\times\wt Y_2$.\bull

\section{Resolving singularities}\label{041010d}
Let $\A$ be a circulant S-ring, $C\in\P_{sgl}(\A)$ a singular class
and $S=S_{min}(C)$. Set
\qtnl{260410i}
\gwr_{\A}(C)=
\css
\gwr_\A(C,\sym(S)), &\text{if $C$ is of composite order,}\\
\gwr_\A(C,\hol(S)),   &\text{otherwise}\\
\ecss
\eqtn
(see Subsection~\ref{161110x}). Statement~(1) of Theorem~\ref{181002b} below
implies that this group is contained in~$\aut(\A)$.

\thrml{181002b}
Let $\A$ be a schurian S-ring over a cyclic group~$G$. Then the group $\aut(\A)$ contains a
$2$-equivalent subgroup $\Gamma\ge G_{right}$ such that for any $C\in\P_{sgl}(\A)$
the following statements hold:
\nmrt
\tm{1} $\gwr_\A(C)\le\Gamma$,
\tm{2} for any $S\in C$ we have $\Gamma^S=\hol(S)$ if $|S|$ is prime, and
 $\Gamma^S=\sym(S)$ otherwise.
\enmrt
\ethrm
\proof Induction on the number $m=m(\A)$ of singular classes of prime order. If $m=0$, then we
are done with $\Gamma=\aut(\A)$. Indeed, statement~(1) immediately follows from statement~(1)
of Theorem~\ref{260410a} because in this case $\rk(\A_S)=2$ and hence $\aut(\A_S)=\sym(S)$
where $S=S_{min}(C)$. But then statement~(2) is a consequence of statement~(1) of
Theorem~\ref{160210c}.\medskip

Suppose $m>0$ and $C$ is a singular class of~$\A$ of prime order. Set
\qtnl{190310a}
\A'=\Ext_C(\A,\Z S)
\eqtn
where $S=S_{min}(C)$. Then $m(\A')=m(\A)-1$ by statement~(2) of Theorem~\ref{200410f}.
So by the induction hypothesis applied to the S-ring $\A'$ the group
$\aut(\A')$ is $2$-equivalent to a group $\Gamma'\ge G_{right}$ such that for any
$C'\in\P_{sgl}(\A')$ statements (1) and (2) hold with $\A$ and $C$ replaced by~$\A'$
and~$C'$ respectively. Set
\qtnl{180111a}
\Gamma=\lg\Gamma',\Delta\rg
\eqtn
where $\Delta=\gwr_\A(C)$ is the group defined in~(\ref{260410i}). Obviously $G_{right}\le\Gamma$. To
complete the proof it suffices to verify that~$\aut(\A)$ is $2$-equivalent to~$\Gamma$
and statements~(1) and~(2) hold. Below we set $L_i=L_i(C)$ and $U_i=U_i(C)$, $i=0,1$.\medskip

To prove that~$\aut(\A)$ is $2$-equivalent to~$\Gamma$ it suffices to verify that the group
$\Gamma_u$ with $u=1_G$, acts transitively on each set $X\in\S(\A)$. However, the S-ring
$\A'$ is schurian by Theorem~\ref{301110a}. Therefore by statement~(3) of Lemma~\ref{011110a} the
transitivity holds when $X\subset U_0$ or $X\subset G\setminus U_1$
because in this case $X\in\S(\A')$ and hence $X\in\orb(\Gamma'_u)$.
Suppose that $X\subset U_1\setminus U_0$. Since $\rk(\A_S)=2$ and $\A'_S=\Z S$,
statements~(2) and~(3) of Lemma~\ref{011110a} imply that there exists a set $Y\in\S(\A)_{U_0}$
such that
$$
X=(L_1\setminus L_0)Y=\bigcup_{x}xL_0Y
$$
where $x$ runs over the representatives of $L_0$-cosets in $L_1$ other than $L_0$. Moreover,
the set $xL_0Y$ is basic in $\A'$ for all~$x$. Therefore the group $\Gamma'_u$ acts transitively
on each of them. On the other hand, the group $\Delta_u$ acts on $S$ as $\hol(S)_u=\aut(S)$,
and hence is transitive on the set $S\setminus\{1\}$. Thus the group $\Gamma_u$ acts transitively
on~$X$ and we are done.\medskip

To prove statements~(1) and (2) let $\wt C\in\P_{sgl}(\A)$. If $\wt C=C$, then statement~(1)
follows because $\gwr_\A(\wt C)=\Delta\le\Gamma$. Otherwise, by statement~(2) of Theorem~\ref{200410f}
we have $\wt C\in\P_{sgl}(\A')$. However for the group~$\Gamma'$ statement~(1) holds
and hence $\gwr_{\A'}(\wt C)\le\Gamma'$. Since $\gwr_{\A'}(\wt C)=\gwr_\A(\wt C)$ and
$\Gamma'\le\Gamma$, statement~(1) is proved.\medskip

Turn to statement~(2).
First suppose that the class $\wt C$ is of composite order. Then $\wt C\ne C$.
By Theorem~\ref{200410f} this implies that $\wt C\in\P_{sgl}(\A')$. Since for the group
$\Gamma'$ statement~(2) holds, for any section $\wt S\in\wt C$ we have
$$
\Gamma^{\wt S}\ge (\Gamma')^{\wt S}=\sym(\wt S)
$$
and we are done. Now suppose that the class $\wt C$ is of prime order. Since by the induction hypothesis
the groups~$\Gamma'$ and $\aut(\A')$ are $2$-equivalent, so are the groups $(\Gamma')^{\wt S}$
and $\aut(\A')^{\wt S}$ where $\wt S\in\wt C$. On the other hand, by the definition of~$\A'$
and Theorem~\ref{261010a} the latter group coincides with $\wt S_{right}$ when $\wt C=C$. Thus
taking into account that for $\Gamma'$ statement~(2) holds, we conclude that
\qtnl{290410a}
(\Gamma')^{\wt S}=\css
\wt S_{right},   &\text{if $\wt C=C$,}\\
\hol(\wt S), &\text{if $\wt C\ne C$.}\\
\ecss
\eqtn
From now on without loss of generality we assume that $\wt S=\wt U_1/\wt U_0$ and
$\wt U_0=\{1_G\}$ where $\wt U_i=U_0(\wt C)$, $i=1,2$: the former due to
Lemma~\ref{160210a} whereas the latter because the natural epimorphism
$\Gamma\to\Gamma^{G/\wt U_0}$ induces an isomorphism from $\Gamma^{\wt S}$ onto
$(\Gamma^{G/\wt U_0})^{\wt S}$. Then $\wt S=\wt U_1$ is a cyclic group of prime order and
$\wt C=\{\wt S\}$.\medskip

To formulate the following auxiliary lemma we need one more notation.
It is easily seen that given a $\wt U_1$-coset $\wt X$ the group
$\hol(\wt U_1)^\gamma\le\sym(\wt X)$
does not depend on the choice of $\gamma\in (G_{right})^{\wt U_1,\wt X}$. We denote this
group by $K(\wt X)$. Clearly, $K(\wt U_1)=\hol(\wt U_1)$.

\lmml{221110a}
Let $K\le\sym(G)$ be the intransitive direct product of the groups $K(\wt X)\le\sym(\wt X)$
where $\wt X$ runs over the set $G/\wt U_1$. Then the group
$\Delta$ normalizes~$K$.
\elmm
\proof 
Since $\wt U_1$ is a block of the group~$\Gamma\ge\Delta$,
it suffices to verify that given $\gamma\in\Delta$ we have
\qtnl{300410g}
K(\wt X)^{\gamma^{\wt X}}=K(\wt X^\gamma)
\eqtn
for all $\wt X\in G/\wt U_1$. To do this suppose first that
$\wt X\not\subset U_1$. Then since $\wt X$ and $U_1$ are blocks of~$\Gamma$, and the number $|\wt X|$
is prime, we see that $|\wt X\cap U_1|\le 1$. However, the group $\Delta=\gwr_\A(C)$
acts trivially outside~$U_1$. Therefore $\gamma$ fixes at least one, and hence all points
of~$\wt X$. Thus in this case equality~(\ref{300410g}) is obvious.\medskip

Let $\wt X\subset U_1$. Then obviously $\wt U_1\subset U_1$. First, suppose that $\wt U_1\subset U_0$. Then $\wt X$ is contained in
some set $X\in U_1/U_0$. By the definition of $\gwr_\A(C)$ (see~\eqref{231110a}) this implies that
$\gamma^{U_1}\in ((U_1)_{right})^{X,X'}$ where $X'=X^\gamma$. This implies
that $\gamma^{\wt X}\in (G_{right})^{\wt X,\wt X'}$ where $\wt X'=\wt X^\gamma$,  whence (\ref{300410g}) follows.
Finally, let $\wt U_1\not\subset U_0$. Since the index of $U_0$ in $U_1$ is prime, this implies
that $\wt U_1U_0=U_1$. On the other hand, since $|\wt U_1|$ is
prime, we have $\wt U_1\cap U_0=\{1\}$. Thus $S=U_1/U_0$ is a multiple of $\wt U_1/1=\wt S$,
and hence $\wt C=C$. This implies that $U_0=L_0=1$, and hence $\wt X=U_1=L_1$.
Thus
$$
K(\wt X)=\hol(U_1)=\hol(L_1)=\gwr_\A(C)^{U_1}=\Delta^{U_1}
$$
and equality~(\ref{300410g}) holds because obviously $(\Delta^{U_1})^{\gamma^{U_1}}=\Delta^{U_1}$.\bull\medskip

To complete the proof we observe that the $\wt U_1$-cosets form an imprimitivity system
for the group~$\Gamma'\ge G_{right}$. From~(\ref{290410a}) it follows that $\Gamma'$
normalizes~$K$. By Lemma~\ref{221110a} and~\eqref{180111a} this implies that the group
$\Gamma$, and hence the group $\Gamma_{\{\wt S\}}$, normalizes~$K$. It follows that
$\Gamma^{\wt S}$ normalizes~$K^{\wt S}$  in $\sym(\wt S)$. Therefore taking into
account~(\ref{290410a}) and the fact that $\Delta^{\wt S}=\hol(\wt S)$ for $\wt C=C$,
we have
\qtnl{290410c}
\hol(\wt S)\le \Gamma^{\wt S}\le N(K^{\wt S})=N(\hol(\wt S))
\eqtn
where $N(\cdot)$ denotes the normalizer in $\sym(\wt S)$.
However, $\wt S$ being a cyclic group of prime order is a characteristic
subgroup of its holomorph. It follows that the group in the right-hand side of~(\ref{290410c})
coincides with $N(\wt S_{right})=\hol(\wt S)$. Thus $\Gamma^{\wt S}=\hol(\wt S)$ and we are done.\bull

\section{Proof of Theorem~\ref{240810a}}\label{051010a}

Below under a {\it section} of a transitive group
$$
\Gamma\le\sym(V)
$$
we mean a permutation group $\Gamma^{X/E}$ where $X$ is a block of~$\Gamma$ and $E$ is a $\Gamma$-invariant
equivalence relation. It is easily seen that this section is transitive; it is primitive whenever
each block of~$\Gamma$ properly contained in~$X$ is contained in a class of~$E$.\medskip

Let now $\Gamma$ be as in Theorem~\ref{240810a}. Then any section of~$\Gamma$ is permutationally
isomorphic to a section $\Gamma^S$ where $S\in\scc(G)$. Therefore statement~(1) immediately
follows from Theorem~\ref{160210c} and the following lemma.

\lmml{270810b}
Any composition factor of a transitive group is isomorphic to a composition factor of some of
its primitive sections.
\elmm
\proof Let $\Gamma\le\sym(V)$ be a transitive group. Take a minimal $\Gamma$-invariant
equivalence relation~$E$. Then by the Jordan-H\"older theorem any composition factor of $\Gamma$ is isomorphic to
either a composition factor of the transitive group $\Gamma^{V/E}$, or a
composition factor of the kernel $\Gamma_E$ of the epimorphism from $\Gamma$ onto $\Gamma^{V/E}$.
Since obviously any primitive section of~$\Gamma^{V/E}$ is permutationally isomorphic
to a primitive section of~$\Gamma$, by induction it suffices to verify that any composition
factor of $\Gamma_E$ is isomorphic to a composition factor of the group $\Gamma^X$ where $X\in V/E$.
However, $\Gamma_E$ is a subdirect product of the groups $(\Gamma_E)^X$, $X\in V/E$.
Therefore any composition factor of $\Gamma_E$ is isomorphic to a composition factor
of a any of the groups $(\Gamma_E)^X$. Since the latter group is normal in $\Gamma^X$,
we are done.\bull\medskip

To prove the sufficiency in statement~(2) suppose that every alternating composition factor
of the group~$\Gamma$ is of prime degree. Denote by $\A$ the S-ring over the group~$G$ associated
with $\Gamma$. Then $\Gamma=\aut(\A)$. We claim that a class $C\in\P_{sgl}(\A)$ of composite degree is in fact
of degree~$4$. To prove this claim let $S\in C$. Then by Theorem~\ref{160210c} we have
$\Gamma^S=\sym(S)$. Besides, if $S=U/L$, then the group $\Gamma_E$ where $E=E_U$, is obviously
normal in~$\Gamma$. Therefore
$$
S_{right}\le(\Gamma_E)^S\,\trianglelefteq\,\Gamma^S=\sym(S).
$$
This implies that either $|S|=4$ or $(\Gamma_E)^S\in\{\alt(S),\,\sym(S)\}$. By the Jordan-H\"older
theorem any composition factor of $(\Gamma_E)^S$ is isomorphic to a composition factor of $\Gamma$.
Thus due to the supposition we conclude that $|S|=4$, and the claim is proved.\medskip

By the claim and Theorem~\ref{181002b} one can find a $2$-equivalent to~$\Gamma$ group $\Gamma'$
such that $G_{right}\le\Gamma'\le\sym(G)$ and $(\Gamma')^S=\hol(S)$ for any section $S$
belonging to a singular class $C\in\P_{sgl}(\A)$ of degree other than~$4$. Since the
group $\sym(4)$ is solvable, we conclude from Theorem~\ref{160210c} that any primitive section
of the group $\Gamma'$ is solvable. By Lemma~\ref{270810b} this implies that the latter group
is solvable. The sufficiency in statement~(2) is proved.\medskip

To prove the necessity suppose that the group $\Gamma$ is $2$-equivalent to a solvable group
$\Gamma'$ containing $G_{right}$. Let $\Gamma$ have a composition factor isomorphic to the group
$\alt(m)$ where $m$ is a composite number. Then obviously $m>4$. By Lemma~\ref{270810b}
there exists a primitive section~$S\in\scc(\A)$ such that the group~$\Gamma^S$ has a composition
factor isomorphic to $\alt(m)$. Since $m>4$, the latter group is not solvable. Therefore
by Theorem~\ref{160210c} we have
$$
\Gamma^S=\sym(S)\quad\text{and}\quad m=|S|.
$$
On the other hand, the groups $\Gamma^S$ and $(\Gamma')^S$ are $2$-equivalent. So the
latter one is primitive. Since it is also solvable, its degree is a prime power~\cite[Theorem~11.5]{W64}. Taking
into account that this group contains the regular cyclic subgroup $S_{right}$, we conclude
by \cite[Theorem~27.3]{W64} that either $|S|$ is prime, or $|S|=4$. However, $m=|S|$
is a composite number greater than~$4$. Contradiction.\bull

\section{Proof of Theorem~\ref{190210a}}

We deduce Theorem~\ref{190210a} from Theorem~\ref{300410w} below the proof of which uses
the following auxiliary lemma.

\lmml{190210b}
Let $G=\prod_{i\in I}G_i$ and $\Delta=\prod_{i\in I}\Delta_i$ where $G_i$
is a cyclic group and $(G_i)_{right}\le \Delta_i\le\hol(G_i)$,
$i\in I$. Then $\Delta'=\Delta$ for any subgroup $\Delta'$ of the group~$\Delta$
that is $2$-equivalent to the latter.
\elmm
\proof By the lemma hypothesis the $2$-orbits
of $\Delta$ and $\Delta'$ are the same. On the other hand, denote by $X_i$ the $2$-orbit
of $\Delta_i$ containing the pair $(1_{G_i},g_i)$ where $g_i$ is a generator of the group~$G_i$.
Then $\Delta_i$ acts regularly on~$X_i$ because $\Delta_i\le\hol(G_i)$.
Therefore the set $X=\prod_iX_i$ is a $2$-orbit of $\Delta$ on which this group acts regularly.
Thus
$$
|\Delta|=|X|\le |\Delta'|.
$$
Since $\Delta'\le\Delta$, this implies that $\Delta'=\Delta$ and we are done.\bull\medskip

Turning to the proof of Theorem~\ref{190210a} set $S=U/L$. Then by the theorem hypothesis
there exists an integer $k\ge 0$ such that
\qtnl{090310b}
\A_S=\bigotimes_{i=0}^k\A_{S_i}
\eqtn
where $S_i$ is an $\A_S$-subgroup, $S=\prod_{i=0}^kS_i$, the S-ring $\A_{S_i}$ is normal
for $i=0$, and is primitive of rank~$2$ and of degree at least~$3$ for $i\ge 1$.
Denote by $I$ (resp. by $J$) the set of all $i\in\{1,\ldots,k\}$ for which the number
$|S_i|$ is prime (resp. composite).

\thrml{300410w}
In the above notations suppose that the S-ring~$\A$ is schurian and
$\Gamma\le\sym(G)$ is the group the existence of which was stated in Theorem~\ref{181002b}.
Then
\qtnl{090310e}
\Gamma^S=\aut(\A_{S_0})\times\prod_{i\in I}\hol(S_i)\times
\prod_{i\in J}\sym(S_i).
\eqtn
\ethrm
\proof Let $j\in J$. Then by Theorem~\ref{160210c} the class $C_j$ of projectively
equivalent sections that contains~$S_j$ is singular. Set $\Delta_j=\gwr_\A(C_j)$.
Then $\Delta_j\le\Gamma$ by the choice of~$\Gamma$ (see statement~(1) of Theorem~\ref{181002b}).
Moreover,
$$
(\Delta_j)^{S_j}=\sym(S_j)\quad\text{and}\quad
(\Delta_j)^{S_i}\le(S_i)_{right},\quad i\ne j
$$
(the equality follows from the definition of the group $\gwr_\A(C_j)$ and
statement~(2) of Theorem~\ref{260410a} whereas the inclusion is the consequence of 
statement~(3) of the same theorem).
Since obviously $\Gamma^S\ge\lg\Delta_j^S,S_{right}\rg$ and $S_{right}$ is the direct
product of the groups $(S_i)_{right}$ over all~$i$, it follows that
$$
\Gamma^S\ge \sym(S_j)\times \{\id_{S_{j'}}\}.
$$
where $S_{j'}$ is the direct product of the groups $S_i$, $i\ne j$. Thus
\qtnl{040510h}
\Gamma^S=\Gamma^{S_{I^*}}\times\prod_{j\in J}\sym(S_j)
\eqtn
where $I^*=I\cup\{0\}$ and $S_{I^*}$ is the product of all $S_i$ with $i\in I^*$.\medskip

By the theorem hypothesis the groups $\Gamma$ and $\aut(\A)$ are $2$-equivalent.
Therefore so are the groups $\Gamma^{S'}$ and $\aut(\A)^{S'}$ for any section $S'\in\scc(\A)$.
Due to~(\ref{090310b}) and the schurity of~$\A$ this implies that
$$
\Gamma^{S_{I^*}}\twoe\aut(\A)^{S_{I^*}}\twoe\aut(\A_{S_{I^*}})=\prod_{i\in I^*}\aut(\A_{S_i})
\twoe \prod_{i\in I^*}\aut(\A)^{S_i}\twoe \prod_{i\in I^*}\Gamma^{S_i}.
$$
Thus by (\ref{040510h}) and Lemma~\ref{190210b} to complete the proof of the theorem it suffices
to verify that given $i\in I^*$ we have
\qtnl{040510l}
\Gamma^{S_i}=\css
\hol(S_i),   &\text{if $i\ne 0$,}\\
\aut(\A_{S_0}),       &\text{if $i=0$.}\\
\ecss
\eqtn
To do this let $i\in I^*$. First suppose that $i=0$. Then the normality of the S-ring $\A_{S_0}$
implies that
$$
\Gamma^{S_0}\le\aut(\A)^{S_0}\le\aut(\A_{S_0})\le\hol(S_0).
$$
On the other hand, we have
$$
\Gamma^{S_0}\twoe\aut(\A)^{S_0}\twoe\aut(\A_{S_0}).
$$
Thus $\Gamma^{S_0}=\aut(\A_{S_0})$ by Lemma~\ref{190210b} for $I=\{0\}$ and
$\Delta=\aut(\A_{S_0})$.\medskip

Suppose that $i>0$. Then $i\in I$ and the number $|S_i|$ is prime. Therefore without loss
of generality we can assume that the section $S_i$ is not singular (see statement~(2)
of Theorem~\ref{181002b}). Then by Theorem~\ref{160210c} we have
$$
\Gamma^{S_i}\le\aut(\A)^{S_i}\le\hol(S_i).
$$
On the other hand, since $\A$ is schurian and $\rk(\A_{S_i})=2$, we obtain that
$$
\Gamma^{S_i}\twoe\aut(\A)^{S_i}\twoe\aut(\A_{S_i})\twoe\hol(S_i).
$$
Thus $\Gamma^{S_i}=\hol(S_i)$ by Lemma~\ref{190210b} for $I=\{i\}$ and
$\Delta=\hol(S_i)$.\bull\medskip

Let us turn to the proof of Theorem~\ref{190210a}. Since the schurity is preserved under
taking restrictions and factors, the necessity is obvious. To prove the
sufficiency suppose that the S-rings $\A_U$ and $\A_{G/L}$ are schurian.
Denote by $\Gamma_0$ and $\Gamma_U$ the groups the existence of which is provided by
Theorem~\ref{181002b} applied to the S-rings $\A_{G/L}$ and $\A_U$ respectively.
Then for the groups $\Delta_0=\Gamma_0$ and $\Delta_U=\Gamma_U$
the first two inclusions in~\eqref{100111a} and condition~\eqref{100111c}
 are  satisfied. By Corollary~\ref{260210a}  to complete the proof it suffices
to verify that the third equality in~\eqref{100111a} also holds for them. However, this
is true due to Theorem~\ref{300410w} applied both to the S-ring $\A_{G/L}$ with $\Gamma=\Delta_0$,
$S=U/L$, and to the S-ring $\A_U$ with $\Gamma=\Delta_U$, $S=U/L$.\bull

\section{Circulant S-rings with $\Omega(n)\le 4$}\label{051010b}

The main goal of the section is to study non-schurian S-rings over a cyclic group of
order~$n$ with $\Omega(n)\le 4$.

\thrml{140510a}
Any S-ring over a cyclic group of order~$n$ with $\Omega(n)\le 3$ is schurian.
\ethrm
\proof The required statement is trivial for $n=1$. Let $\Omega(n)>0$ and
$\A$ an S-ring over a cyclic group $G$ of order~$n$. If the radical of~$\A$ is trivial,
then we are done by Corollary~\ref{170510b}. Let now the radical of $\A$ is not trivial. Then by statement~(1) of Theorem~\ref{mbnulc}
there exists an $\A$-section $S=U/L$ of the group $G$ such that the S-ring
$\A$ satisfies the $S$-condition nontrivially. By the induction hypothesis the S-rings $\A_U$
and $\A_{G/L}$ are schurian. Moreover, since $\Omega(n)\le 3$, we have $\Omega(|S|)\le 1$.
Therefore $\rad(\A_S)=1$. Thus $\A$ is schurian by Theorem~\ref{190210a}.\bull\medskip

The following auxiliary lemma will be used in the analysis of the four primes case.
We fix a cyclic group $G$ of order $pqr$ where $p$, $q$ and $r$ are primes.
The subgroups of $G$ of orders $pq$ and $r$ are denoted by $M$ and $N$ respectively.
We also fix a section $S$ of $G$ of order~$pq$; obviously, $S=M/1$ or $S=G/N$.

\lmml{071210a}
Let $\A$ be an S-ring over the group $G$ that is not a proper wreath product.
Suppose that $S\in\scc(\A)$, $|S|\ne 4$ and $\A_S$ is the wreath product of the S-ring over the
group of order~$p$ by the S-ring over the group of order~$q$. Then the set of orders
of proper $\A$-subgroups equals~$\{p,r,pq,pr\}$ and
\nmrt
\tm{1} if $S=M/1$, then $q\ne r$,
\tm{2} if $S=G/N$, then $p\ne r$.
\enmrt
Moreover, either $\A$ is normal and then $p=q$, or $\A$ satisfies the $M'/N'$-condition
where $M'$ and $N'$ are the subgroups of~$G$ of orders $pr$ and $p$ respectively.
\elmm
\proof First, we find the set $\H(\A)$ and prove that statements~(1) and~(2) hold
assuming $G=M\times N$ with $M,N\in\H(\A)$. Indeed, in this case obviously  $r\not\in\{p,q\}$
and the section $S$ is projectively equivalent to both $M/1$ and $G/N$. Therefore
both statements~(1) and~(2) hold. Next, since by Theorem~\ref{261010a} the S-ring
$\A_M$ is Cayley isomorphic to $\A_S$, the hypothesis of the lemma on the latter
implies that $\H(\A_M)=\{1,N',M\}$. Besides, $\H(\A)=\H(\A_M)\H(\A_N)$ by Lemma~\ref{130209d}.
Thus $\H(\A)=\{1,N,N',M',M,G\}$ and we are done.\medskip


Suppose that the S-ring~$\A$ is normal. Then it is dense by
statement~(1) of Theorem~\ref{030910a}. This implies that $p=q$ because otherwise
$|\H(\A_S)|=3$ which contradicts the lemma hypothesis. Thus if $G$ is not a $p$-group,
then $G=M\times N$ with $M,N\in\H(\A)$, and we are done by the previous paragraph.
However, if $G$ is a $p$-group, then $p\ne 2$ because by the hypothesis $|S|\ne 4$. Then by Corollary~\ref{091210v}
this implies that the S-ring $\A_S$ is normal, which is impossible due to
statement~(3) of Theorem~\ref{030910a}.\medskip

Finally, suppose that the S-ring $\A$ is not normal. If $\rad(\A)=1$, then by statement~(2) of
Theorem~\ref{mbnulc} this implies that
\qtnl{101210a}
\A=\A_{\wt M}\otimes\A_{\wt N}
\eqtn
for some $\A$-groups $\wt M$ and $\wt N$ with $\Omega(|\wt M|)=2$ and $\Omega(|\wt N|)=1$ (here
we use the fact that $\rk(\A)\ge 3$ by the lemma hypothesis). However, if $|\wt N|\in\{p,q\}$,
then it is easily seen that $\A_S$ is a nontrivial tensor product which is impossible. Thus
$\wt M=M$ and $\wt N=N$. Therefore the section $M/1$ is projectively equivalent to~$S$, and so
$N'\in\H(\A)$ and $\A_M=\A_{N'}\wr\A_{M/N'}$. Thus
due to~\eqref{101210a} the S-ring $\A$ satisfies the $M'/N'$-condition which contradicts
the assumption $\rad(\A)=1$ by statement~(1) of Theorem~\ref{mbnulc}.\medskip

Let $\rad(\A)\ne 1$.
Then by statement~(1) of Theorem~\ref{mbnulc} there exist $\A$-groups $\wt M$ and $\wt N$
such that $1<\wt N\le\wt M<G$ and $\A$ satisfies the $\wt M/\wt N$-condition.
Since $\A$ is not a proper wreath product, we also have $\wt M\ne\wt N$. By the same
reason both $\A_{\wt M}$ and $\A_{G/\wt N}$ are not proper wreath products (because
otherwise $\A=\A_{\wt N}\wr\A_{G/\wt N}$ in the former case, and $\A=\A_{\wt M}\wr\A_{G/\wt M}$
in the latter case).\medskip

Let us assume that $S=M/1$ (the case $S=G/N$ can be treated in a
similar way). Then obviously $M$ is an $\A$-group. Moreover,
\qtnl{201210a}
\wt N\subset M.
\eqtn
Indeed, otherwise $\wt N\cap M=1$, and hence $G=M\times\wt N$. But then $\A_{G/\wt N}$ being
a Cayley isomorphic to~$\A_M$, is a proper wreath product, which is impossible by above. Next,
since $\A_{\wt M}$ is also not a proper wreath product, we have
\qtnl{201210b}
\wt M\ne M.
\eqtn
Now, taking into account that $\H(\A_M)=\{1,N',M\}$, we deduce from~\eqref{201210a} that
$\wt N=N'$. Therefore $\wt M\in\{M,M'\}$ and due to~\eqref{201210b} we have
$\wt M=M'$, and hence $M'\ne M$. The latter implies that $q\ne r$. Finally,
to find the set $\H(\A)$ we observe that since the S-ring $\A_{M'}=\A_{\wt M}$ is of rank
at least~$3$ and not a proper wreath product, Theorem~\ref{mbnulc} implies that
$\A_{M'}$ is either normal or a nontrivial tensor product. So it is dense
(in the former case this follows from statement~(1) of Theorem~\ref{030910a}),
and hence $\H(\A_{M'})\supset\{N',N\}$. Therefore
$\H(\A)$ is as required: if $N=N'$, then this is obvious, whereas if $N\ne N'$, then
$G=M\times N$ and this is true by the statement in the first paragraph.\bull

\thrml{211210a}
Let $\A$ be a non-schurian S-ring over a cyclic group $G$ of order~$n$ with $\Omega(n)=4$.
Then $\A$ satisfies nontrivially some $S$-condition. Moreover,
for any such $\A$-section $S=U/L$ we have
\nmrt
\tm{1} $|L|$ and $|G/U|$ are primes, and $|S|\ne 4$,
\tm{2} $\A_S$ is a proper wreath product,
\tm{3} $\A_U$ and $\A_{G/L}$ are not proper wreath products and cannot be normal
simultaneously.
\enmrt
\ethrm
\proof By Lemma~\ref{170510b} the radical of the S-ring~$\A$ is not trivial. Then by
statement~(1) of Theorem~\ref{mbnulc} there exists an $\A$-section $S=U/L$ of the group~$G$ such
that the S-ring $\A$ satisfies the $S$-condition nontrivially. It follows that
$|L|$ and $|G/U|$ are primes. Therefore by Theorem~\ref{140510a} the S-rings $\A_U$
and $\A_{G/L}$ are schurian. By Theorem~\ref{190210a} this implies that the S-ring
$\A_S$ can be neither normal nor of rank~$2$. Since any S-ring over a cyclic group of order~$4$
is obviously of rank~$2$ or normal, statement~(1) follows. Moreover, by Theorem~\ref{mbnulc}
the S-ring $\A_S$ is a proper generalized wreath product, and hence a proper wreath
product. Statement~(2) is proved.\medskip

Next, the schurity of the S-rings $\A_U$ and $\A_{G/L}$ imply that the groups
$\aut(\A_U)^S$ and $\aut(\A_{G/L})^S$ are $2$-equivalent. If in addition both of these
S-rings are normal, then the above groups are between $S_{right}$ and~$\hol(S)$. Thus they
are equal. By Corollary~\ref{260210a} this implies that $\A$ is schurian. Contradiction.
This proves the second part of statement~(3). To prove the first one suppose that
$$
\A_U=\A_H\wr\A_{U/H}
$$
for some proper $\A$-subgroup $H$ of~$U$ (the case when $\A_{G/L}$ is a proper wreath
product can be treated in a similar way). Then any $\A$-subgroup of $U$ either contains $H$
or is a subgroup of~$H$. Since $|L|$ is prime, this implies that $L\le H$. However, then
$\A$ satisfies the $H/L$-condition which is impossible by statement~(1) of
this theorem.\bull\medskip

Let us fix a cyclic group $G$ of order $n=p_1p_2p_3p_4$ with $p_i$'s being primes. Then for
any divisor~$m$ of $n$ there is a unique subgroup of $G$ of order~$m$. Denote by
$U$, $V$, $G_1$, $H$, $G_2$, $K$, $L$  subgroups of $G$ such that
$$
|U|=p_1p_2p_3,\quad |V|=p_1p_3p_4,
$$
$$
|G_1|=p_1p_2,\quad |H|=p_1p_3,\quad |G_2|=p_3p_4,
$$
$$
|K|=p_1,\quad |L|=p_3.
$$
It is easily seen that if $p_1\ne p_3$ and $p_2\ne p_4$, then $K\ne L$, $U\ne V$ and
these subgroups form a sublattice of the lattice of all subgroups of $G$, as in Fig.1 below
(if $p_2=p_3$, then $G_1=H$, whereas if $p_1=p_4$, then $G_2=H$):
\begin{figure}[h]
\vspace{5mm}
\unitlength 1.00mm
\linethickness{0.4pt}
\begin{picture}(111.00,152.00)
\put(89.00,150.00){\circle*{2.67}}
\put(99.00,140.00){\circle*{2.67}}
\put(79.00,140.00){\circle*{2.67}}
\put(69.00,130.00){\circle*{2.67}}
\put(89.00,130.00){\circle*{2.67}}
\put(109.00,130.00){\circle*{2.67}}
\put(99.00,120.00){\circle*{2.67}}
\put(79.00,120.00){\circle*{2.67}}
\put(89.00,110.00){\circle*{2.67}}
\multiput(89.00,110.00)(-0.12,0.12){84}{\line(0,1){0.12}}
\multiput(79.00,120.00)(-0.12,0.12){84}{\line(0,1){0.12}}
\multiput(69.00,130.00)(0.12,0.12){84}{\line(0,1){0.12}}
\multiput(79.00,140.00)(0.12,0.12){84}{\line(0,1){0.12}}
\multiput(89.00,150.00)(0.12,-0.12){84}{\line(0,-1){0.12}}
\multiput(99.00,140.00)(0.12,-0.12){84}{\line(0,-1){0.12}}
\multiput(109.00,130.00)(-0.12,-0.12){84}{\line(0,-1){0.12}}
\multiput(99.00,120.00)(-0.12,-0.12){84}{\line(0,-1){0.12}}
\multiput(79.00,120.00)(0.12,0.12){84}{\line(0,1){0.12}}
\multiput(89.00,130.00)(-0.12,0.12){84}{\line(0,1){0.12}}
\multiput(99.00,140.00)(-0.12,-0.12){84}{\line(0,-1){0.12}}
\multiput(89.00,130.00)(0.12,-0.12){84}{\line(0,-1){0.12}}
\put(89.00,108.00){\makebox(0,0)[ct]{1}}
\put(89.00,152.00){\makebox(0,0)[cb]{$G$}}
\put(77.00,140.00){\makebox(0,0)[rc]{$U$}}
\put(91.00,130.00){\makebox(0,0)[lc]{$H$}}
\put(77.00,120.00){\makebox(0,0)[rc]{$K$}}
\put(101.00,140.00){\makebox(0,0)[lc]{$V$}}
\put(111.00,130.00){\makebox(0,0)[lc]{$G_2$}}
\put(101.00,120.00){\makebox(0,0)[lc]{$L$}}
\put(67.00,130.00){\makebox(0,0)[rc]{$G_1$}}
\end{picture}
\vspace{-105mm}
\caption{}
\label{lf}
\end{figure}




\thrml{051010c}
Let $\A$ be a non-schurian S-ring over the group $G$ as above. Suppose that $\A$ satisfies
the $U/L$-condition. Then $|U/L|\ne 4$ and
\nmrt
\tm{1} $p_1\ne p_3$ and $p_2\ne p_4$ and the lattice of $\A$-subgroups of~$G$ is as
in Fig.1,
\tm{2} $\A_{U/L}=\A_{H/L}\wr\A_{U/H}$,
\tm{3} the S-rings $\A_U$ and $\A_{G/L}$ cannot be normal simultaneously; moreover,
if one of them is normal, then $p_1=p_2$,
\tm{4} if $\A_U$ is not normal, then $\A_U=\A_H\wr_{H/K}\A_{U/K}$, whereas if
$\A_{G/L}$ is not normal, then $\A_{G/L}=\A_{V/L}\wr_{V/H}\A_{G/H}$,
\tm{5} if both  $\A_U$ and $\A_{G/L}$ are not normal, then both $\A_V$ and $\A_{G/K}$
are not normal and $\A=\A_V\wr_{V/K}\A_{G/K}$,
\enmrt
\ethrm
\proof Let $S=U/L$. Then the inequality $|S|\ne 4$, statement~(2) and the first part of statement~(3)
hold by Theorem~\ref{211210a}. Besides, by the same theorem the hypothesis of
Lemma~\ref{071210a} is satisfied for the section $S$ and both the S-ring $\A_U$
and $(p,q,r)=(p_1,p_2,p_3)$, and the S-ring $\A_{G/L}$ and $(p,q,r)=(p_1,p_2,p_4)$.
Thus by this lemma all the rest except for statement~(5) follow. To prove the latter suppose that
the S-rings $\A_U$ and $\A_{G/L}$ are not normal. Then by statement~(4) we have
\qtnl{270111x}
\A_U=\A_H\wr_{H/K}\A_{U/K}\quad\text{and}\quad
\A_{G/L}=\A_{V/L}\wr_{V/H}\A_{G/H}.
\eqtn
Take $X\in\S(\A)_{G\setminus V}$. If $X\subset U\setminus V$, then obviously
$X\subset U\setminus H$. By the left-hand side equality in~\eqref{270111x} this implies
that $K\le\rad(X)$. Let now $X\subset G\setminus U$. Then $L\le\rad(X)$ because $\A$
satisfies the $U/L$-condition by the theorem hypothesis. However,
$\pi_L(X)\subset \pi_L(G)\setminus \pi_L(V)$. By the right-hand side equality in~\eqref{270111x}
this implies that $\pi_L(H)\le \rad(\pi_L(X))$. Therefore $K\le H\le\rad(X)$.
Thus
$$
\A=\A_V\wr_{V/K}\A_{G/K}.
$$
Next, by statement~(1) of Theorem~\ref{211210a} we have $|V/K|\ne 4$. Besides, the S-ring $\A_V$
satisfies the $H/L$-condition because $\A$ satisfies the $U/L$-condition, and $U\cap V=H$ due
to the inequality $p_2\ne p_4$ proved above. Therefore by statement~(3) of
Theorem~\ref{030910a} the S-ring $\A_V$ can be normal only if $|V/H|=|L|=2$.
But this is not the case because
$$
4\ne |V/K|=p_3p_4=|L|\,|V/H|.
$$
Thus the S-ring $\A_V$ is not normal. The fact that the S-ring $\A_{G/K}$ is not normal
is proved in a similar way.\bull

\section{Example}\label{051010g}
Throughout this section $\Z_n$ is the additive group of integers modulo a positive integer~$n$.\medskip

For any divisor $m$ of~$n$ denote by $i_{m,n}:\Z_m\to\Z_n$ and $\pi_{n,m}:\Z_n\to\Z_m$ the
group homomorphisms taking~$1$ to~$n/m$ and to~$1$ respectively. Using them we identify
with~$\Z_m$ both the subgroup $i_{m,n}(\Z_m)$ and the factorgroup $\Z_n/\ker(\pi_{n,m})$
of~$\Z_n$. Thus every section of $\Z_n$ of order~$m$ is identified with the group~$\Z_m$.
Moreover, the automorphism of~$\Z_n$ afforded by multiplication by~$d$ induces the
automorphism of that section afforded  by multiplication by the same number~$d$.\medskip

If~$\A$ is an S-ring over~$\Z_n$ and $\Z_m$ belongs to~$\S^*(\A)$, then $\A_{\Z_m}$ will
be denoted by $\A_m$ and $\A_{\Z_n/\Z_m}$ by~$\A^{n/m}$.
Let finally $\A_l$ be an S-ring over $\Z_{n_l}$
($l=1,2$) and $(\A_1)^m=(\A_2)_m$ for some~$m$ dividing both~$n_1$ and $n_2$.
Then the unique S-ring~$\A$ over $\Z_{n_1n_2/m}$ from Theorem~\ref{160710a}
is denoted by $\A_1\wr_m\A_2$. We omit~$m$ if~$m=1$. Given a group~$K\le\aut(\Z_n)$ and a
prime~$p$ dividing $n$ we write $K_p$ for the $p$-projection of~$K$ in the sense
of decomposition~\eqref{060910a}, and $\cyc(K,n)$ instead of $\cyc(K,\Z_n)$.\medskip

In what follows we will construct a family of non-schurian S-rings over a cyclic group of order
$n=p_1p_2p_3p_4$ where $p_i$ is a prime and
$$
p_1=p_2=p\ne p_4\quad\text{and}\quad p\ \ \text{divides}\ \ p_3-1.
$$
Let us fix a positive integer~$d$. Suppose that $d\,|\,p-1$ and $M\le\aut(\Z_{p^2p_3})$ a cyclic
group of order~$pd$ such that $|M_p|=pd$ and $|M_{p_3}|=p$. Set
$$
\A_1=\cyc(M,p^2p_3).
$$
We claim that
\qtnl{270111g}
(\A_1)^{p^2}=\cyc(d,p)\wr\cyc(d,p)
\eqtn
where $\cyc(d,p)=\cyc(K,p)$ with $K$ being the subgroup of $\aut(\Z_p)$ of order~$d$.
Indeed, it is easily seen that $(\A_1)^{p^2}=\cyc(M_p,p^2)$. On the other hand, since
$p$ divides $|M_p|$, the group $\rad((\A_1)^{p^2})$ is nontrivial. By statement~(1)
of Theorem~\ref{mbnulc} this implies that $(\A_1)^{p^2}$ is a nontrivial generalized wreath
product. So this S-ring is the wreath product of $(\A_1)^p$ by $(\A_1)_p$.
Since obviously $(\A_1)^p=(\A_1)_p=\cyc(d,p)$, the claim is proved.\medskip

Next, suppose that $d\,|\,p_4-1$. Given two cyclic groups
$M_1,M_2\le\aut(\Z_{pp_4})$ of order~$d$ such that $|(M_i)_p|=|(M_i)_{p_4}|=d$,
$i=1,2$, set
$$
\A_2=\cyc(M_1,pp_4)\wr_{p_4}\cyc(M_2,pp_4).
$$
Then $\rad(\A_2)$ contains the subgroup of order~$p$. Besides, since $p\ne p_4$, Lemma~\ref{130209d}
implies that $\A_2\ge(\A_2)_{p^2}\otimes(\A_2)_{p_4}$. Thus the group $\rad((\A_2)_{p^2})$
is nontrivial. By statement~(1) of Theorem~\ref{mbnulc} this implies that
$(\A_2)_{p^2}$ is a nontrivial generalized wreath product. So
this S-ring is the wreath product of $(\A_2)^p$ by $(\A_2)_p$.
Since $(\A_2)^p=(\A_2)_p=\cyc(d,p)$, we conclude that
\qtnl{270111h}
(\A_2)_{p^2}=\cyc(d,p)\wr\cyc(d,p).
\eqtn

From \eqref{270111g} and \eqref{270111h} it follows that $(\A_1)^{p^2}=(\A_2)_{p^2}$ and we set
\qtnl{210910v}
\A=\A_1\wr_{p^2}\A_2.
\eqtn
We observe that the S-rings $\A_1$, $\cyc(M_1,pp_4)$ and $\cyc(M_2,pp_4)$ are cyclotomic,
and hence dense. The density of the two latter S-rings implies that the S-ring $\A_2$ is also
dense. Thus the set
$\H(\A)$ contains the subgroups $K$, $L$, $G_1$, $H$, $G_2$, $U$ and $V$, of the group
$G=\Z_n$ of orders $p$, $p_3$, $p^2$, $pp_3$, $p_3p_4$, $p^2p_3$ and $pp_3p_4$ (see
Fig.~\ref{lf}).

\thrml{250510x}
The S-ring $\A$ is non-schurian whenever $M_1\ne M_2$.
\ethrm
\proof Suppose that $M_1\ne M_2$. It is easily seen that each of the S-rings $\A_1$,
$\cyc(M_1,pp_4)$ and $\cyc(M_2,pp_4)$ has trivial radical, and is not a proper tensor product
(the latter is true because otherwise $d=1$ and hence $M_1=M_2$). By Theorem~\ref{mbnulc}
this implies that all these S-rings are normal. To prove that the S-ring~$\A$ is not schurian
it suffices to find an element
\qtnl{200910a}
\gamma\in(\Gamma^{U/L})_u\setminus (\Delta^{U/L})_u
\eqtn
where $\Gamma=\aut(\A_1)$, $\Delta=\aut(\A_2)$ and $u=0$ (in the group $U/L$).
Indeed, in this case the group $\Gamma'=\Gamma^{U/L}\cap \Delta^{U/L}$ is a proper subgroup
of $\Gamma^{U/L}$. Besides by the normality of the S-ring~$\A_1$ we have
\qtnl{200910g}
\Gamma^{U/L}\le\hol(U/L).
\eqtn
Thus the groups $\Gamma'$ and $\Gamma^{U/L}$ are not $2$-equivalent by Lemma~\ref{190210b}.
Since obviously $\Gamma^{U/L}$ is $2$-equivalent to $\aut(\A_{U/L})$,
the S-ring~$\A$ is not schurian by Corollary~\ref{250510a}.\medskip

To find $\gamma$ as in~\eqref{200910a} we note that from the definition it follows that $M_1$
is a subgroup of the group
$$
\aut(V/L)=\aut(H/L)\times\aut(G_2/L)=\aut(\Z_p)\times\aut(\Z_{p_4}).
$$
Moreover, if $Q\le\aut(\Z_p)$ and $R\le\aut(\Z_{p_4})$ are the groups of order~$d$, then
\qtnl{200910b}
M_1=\{(x,y)\in Q\times R:\ \varphi_1(x)=y\}
\eqtn
where $\varphi_1$ is an isomorphism from~$Q$ to~$R$. Similarly, $M_2$
is a subgroup of the group
$$
\aut(G/H)=\aut(U/H)\times\aut(V/H)=\aut(\Z_p)\times\aut(\Z_{p_4}),
$$
and
\qtnl{200910c}
M_2=\{(x,y)\in Q\times R:\ \varphi_2(x)=y\}
\eqtn
where $\varphi_2$ is an isomorphism from~$Q$ to~$R$. Let now $\delta\in\Delta_u$. Then
taking into account that the S-rings $(\A_2)_U=\cyc(M_1,pp_4)$ and
$(\A_2)_{G/H}=\cyc(M_2,pp_4)$ are normal (see above), we conclude that
$(\delta)^{V/L}\in M_1$ and $(\delta)^{G/H}\in M_2$. Due to~\eqref{200910b}
and~\eqref{200910c} this implies that
\qtnl{200910d}
\varphi_1(\delta^{H/L})=\delta^{V/H}\quad\text{and}\quad
\varphi_2(\delta^{U/H})=\delta^{V/H}.
\eqtn
Next, since $M_1\ne M_2$, from~\eqref{200910b}
and~\eqref{200910c} it follows that there exist distinct elements $x_1,x_2\in Q$ such
that
\qtnl{200910e}
\varphi_1(x_1)=\varphi_2(x_2).
\eqtn
On the other hand, due to inclusion~\eqref{200910g} one can find an element $\gamma$
in the group $\Gamma_u=M$ for which $\gamma^{H/L}=\gamma^{U/H}=x_1$.
To complete the proof of the theorem let us check that
$\gamma$ satisfies~\eqref{200910a}. Suppose on the contrary that $\gamma^{U/L}\in(\Delta_u)^{U/L}$.
Then there exists $\delta\in\Delta_u$ such that $\delta^{U/L}=\gamma^{U/L}$.
So from~\eqref{200910d}, \eqref{200910e} and the definition of~$\gamma$ it follows that
$$
\varphi_2(\delta^{U/H})=\varphi_1(\delta^{H/L})=\varphi_1(\gamma^{H/L})=\varphi_1(x_1)=\varphi_2(x_2).
$$
This implies that $x_2=\delta^{U/H}=\gamma^{U/H}=x_1$, which is impossible because
$x_1\ne x_2$.\bull

\thrm
Let $n=p_1p_2p_3p_4$ where $p_i$ is an odd prime. Suppose that
$\{p_1,p_2\}\cap\{p_3,p_4\}=\emptyset$ and that $p_1-1$ and $p_4-1$
have a common divisor $d\ge 3$. Then a cyclic group of order $n$ is non-schurian
whenever either $d$ divides both $p_2-1$ and $p_3-1$, or $p_1=p_2$ and $p_1$ divides $p_3-1$.
\ethrm
\proof If the number $d$ divides both $p_2-1$ and $p_3-1$, then $d$ divides $p_i-1$ for
$i=1,2,3,4$. Due to \cite{EP01ae} this implies that there is a non-schurian S-ring over the
group $\Z_n$. Thus this group is non-schurian. To complete the proof suppose that
$p_1=p_2=p$ and $p$ divides $p_3-1$. Using~\eqref{060910a} set
$$
M=\{(x,y)\in\aut(\Z_{p^2p_3}):\ \varphi(x)=y\}
$$
and for $i=1,2$
$$
M_i=\{(x,y)\in\aut(\Z_{pp_4}):\ \varphi_i(x)=y\}
$$
where $\varphi$ is an epimorphism from the subgroup of $\aut(\Z_{p_3})$ of order $pd$
to the subgroup of $\aut(\Z_{p^2})$ of order $p$, and $\varphi_i$ is an isomorphism from
the subgroup of $\aut(\Z_p)$ of order $d$ to the subgroup of $\aut(\Z_{p_4})$ of order $d$.
Then obviously $M$ is a cyclic group of order~$pd$ such that $|M_p|=pd$ and $|M_{p_3}|=p$,
and $M_i$ is a cyclic group of order~$d$ such that $|(M_i)_p|=|(M_i)_{p_4}|=d$.
Since $d\ge 3$, the isomorphisms $\varphi_1$ and $\varphi_2$, and hence
the groups $M_1$ and $M_2$, can be chosen to be distinct. Thus by Theorem~\ref{250510x}
the S-ring $\A$ is non-schurian and we are done.\bull\medskip

The minimal example is $n=5\cdot 5\cdot 11\cdot 13$.

\end{document}